\documentclass[11pt]{article}
\usepackage{amsmath} \usepackage{amssymb} \usepackage{latexsym} \usepackage{enumitem} \usepackage{mathrsfs} \usepackage{comment}
\usepackage{color}
\usepackage{cite}

\usepackage[colorlinks=true,urlcolor=blue,
citecolor=red,linkcolor=blue,linktocpage,pdfpagelabels,bookmarksnumbered,bookmarksopen]{hyperref}

\newtheorem{assumption}{Assumption}

\def\qed{ \ \vrule width.2cm height.2cm depth0cm\smallskip}
\newenvironment{proof}{\noindent {\bf Proof.\/}}{$\qed$\vskip 0.1in}

\newcommand{\ba}{\begin{array}}
\newcommand{\ea}{\end{array}}
\newcommand{\be}{\begin{equation}}
\newcommand{\ee}{\end{equation}}
\newcommand{\bea}{\begin{eqnarray}}
\newcommand{\eea}{\end{eqnarray}}
\newcommand{\beaa}{\begin{eqnarray*}}
\newcommand{\eeaa}{\end{eqnarray*}}

\def\dbE{\mathbb{E}}
\def\dbF{\mathbb{F}}
\def\dbG{\mathbb{G}}

\def\dbL{\mathbb{L}}

\def\dbP{\mathbb{P}}
\def\dbR{\mathbb{R}}

\def\rme{\mathrm{e}}
\def\rmd{\mathrm{d}}
\def\mcu{\mathcal{U} }
\def\mcv{\mathcal{V} }

\def\si{\sigma}

\def\Th{\Theta}

\def\O{\Omega}
\def\cB{{\cal B}}

\def\cF{{\cal F}}
\def\cG{{\cal G}}
\def\cH{{\cal H}}

\def\cK{{\cal K}}
\def\cL{{\cal L}}

\def\cN{{\cal N}}

\def\cP{{\cal P}}

\def\cU{{\cal U}}

\def\cW{{\cal W}}

\def\no{\noindent}

\def\q{\quad}

\def\pa{\partial}

\def\qed{ \hfill \vrule width.25cm height.25cm depth0cm\smallskip}

\newcommand{\basa}{\begin{assumption}}
\newcommand{\easa}{\end{assumption}}

\newcommand{\bas}{\begin{assum}}
\newcommand{\eas}{\end{assum}}

\def\liminf{\mathop{\underline{\rm lim}}}

\def\pa{\partial}

\def\1{{\bf 1}}

\def\:{\!:\!}
\def\reff{\eqref}
\def \proof{{\noindent \bf Proof.\quad}}

\def\Oin{\Omega_{\text{input}}}
\def\Oout{\Omega_{\text{output}}}
\def\Ocan{\Omega_{\text{canon}}}
\def\Oto{\Omega_{\text{total}}}
\def\Qin{Q_{\text{input}}}
\def\oin{\omega_{\text{input}}}

at 9pt

\definecolor{alp}{rgb}{0.0, 0.5, 0.0}

\newtheorem{thm}{Theorem}[section]
\newtheorem{lem}[thm]{Lemma}

\newtheorem{rem}[thm]{Remark}
\newtheorem{eg}[thm]{Example}
\newtheorem{defn}[thm]{Definition}
\newtheorem{assum}[thm]{Assumption}

\begin{document}

\title{\bf On Discounted Infinite-Time Mean Field Games  } 
\author{ Yongsheng Song$^*$ and Zeyu Yang\thanks{State Key Laboratory of Mathematical Sciences, Academy of Mathematics and Systems Science, Chinese Academy of Sciences, Beijing 100190, China, and School of Mathematical Sciences, University of Chinese Academy of Sciences, Beijing 100049, China. E-mails: yssong@amss.ac.cn (Y. Song), yangzeyu@amss.ac.cn (Z. Yang).}  }
\date{\today}

\maketitle

\begin{abstract} 
  In this paper, we study the infinite-time mean field games with discounting, 
  establishing an equilibrium where individual optimal strategies collectively regenerate the mean-field distribution.
  To solve this problem, we partition all agents into a representative 
  player and the social equilibrium.
   When the optimal strategy of the representative 
  player has the same feedback form as the strategy in the social equilibrium, we say that the system achieves a Nash 
  equilibrium.
  We construct a Nash equilibrium using the 
  stochastic maximum principle and infinite-time forward-backward stochastic differential equations (FBSDEs). 
  By employing  elliptic master equations, a class of distribution-dependent elliptic partial differential equations (PDEs), 
  we provide a representation for the Nash equilibrium strategies.
  We prove the Yamada-Watanabe type theorem and show  weak uniqueness for infinite-time FBSDEs.
  Furthermore, we prove that the solutions to a system of infinite-time FBSDEs
  can be employed to construct viscosity solutions for a class of distribution-dependent elliptic PDEs.
\end{abstract}

\no{\bf Keywords.}  discounted infinite-time mean field games, infinite-time forward-backward equations, weak uniqueness, elliptic 
master equations.

\vfill\eject


\section{Introduction}
\label{sect-Introduction}
\setcounter{equation}{0} 

The study of mean field games was initiated independently by Lasry-Lions
\cite{lasry2006jeux,lasry2006mean,lasry2007mean} and Huang-Malhamé-Caines \cite{huang2006large} as an analysis
of limit models for symmetric weakly interacting $(N+1)-$player differential games.
It is noteworthy that current theoretical frameworks are primarily developed for finite-time problems, while infinite-time scenarios remain considerably underdeveloped.
We refer the
reader to \cite{carmona2018probabilistic} for a comprehensive exposition on this subject.

In this paper, we consider a generalized framework for mean field games, which extends classical finite-time settings to
discounted infinite-time mean field games.
In our framework, a representative player interacts with a continuum of other players 
(also referred to as the population or social equilibrium).
Let $\mu_0\in \cP(\dbR)$ be the initial state distribution in the society,
and let $(\mu_t)_{t\geq 0}$ denote the population distribution flow starting at $\mu_0$.
For the player($x$), a representative player with initial state $x$,
the dynamic of its state  is given by
\begin{equation}
    \label{eq: constraint}
  X^{x,\beta}_t=x+\int_0^tb(X^{x,\beta}_t,\mu_s,\beta_s)\rmd s+B_t,
\end{equation}
where  $\beta$, an $\dbR$-valued progressively measurable stochastic process,
is the strategy of the player($x$).
Denote by $\beta^x$ the optimal strategy of the representative player($x$) 
derived from the  stochastic optimization problem:
\begin{equation}
  \min_{\beta} J^\mu(\beta)\triangleq\mathbb{E}\left[\int_0^\infty e^{-rt}f(X_t^{x,\beta},\mu_t,\beta_t)\rmd t\right].
\end{equation}
Since the population consists of a multitude of individuals, the macroscopic distribution should satisfy
\begin{equation}
  \label{eq: relationship}
  \mu_t(\cdot)=\int \dbP(X_t^{x,\beta^x}\in \cdot)\mu_0(\rmd x).
\end{equation}
This is exactly a fixed-point problem.

Note that our definition of mean field games differs from that in \cite{bayraktar2023solvability} and 
from those of earlier finite-time mean field games. 
In \cite{bayraktar2023solvability}, 
all agents are identical and thus represented by a single 
representative player. For a given measure flow $(\mu_t)_{t\ge 0}$, the representative player seeks to minimize
\begin{equation}
  J^\mu(\alpha)\triangleq\mathbb{E}\left[\int_0^\infty e^{-rt}f(t,X_t,\mu_t,\alpha_t)\rmd t\right]
\end{equation}
under the constraint
\begin{equation}
  \begin{cases}
  \rmd X_t=b\left(t,X_t,\mu_t,{\alpha}_t\right)\rmd t+\sigma \rmd B_t, \\
 X_0=\xi. 
\end{cases}
\end{equation}
Then, we require that the law of $X_t$ coincides with $\mu_t$, which means that we need to find a fixed point.

In the conventional definition of mean field games, the primary focus on the distribution flow $(\mu_t)_{t\geq 0}$
stems from its role as the limit of the Nash equilibrium of an $N$-player game. 
Within this framework, the master equation also serves as a crucial tool: see \cite{10.1214/19-EJP298} and \cite{10.1214/19-AOP1359}. 
In the study by Buckdahn et al. \cite{buckdahn2017mean}, the individual state process and the  population state process are 
treated separately, allowing for an investigation into various properties of the parabolic master equation. 
Inspired by their approach, our definition of mean field games explicitly decouples the representative player 
from the population. This separation offers two key advantages: first, it provides a mathematically  
interpretation of the game problem even in the absence of a realistic $N$-player game context; second, 
the value function $V(x,\mu_0)$ of the representative player($x$) 
directly embodies the relationship between the individual 
state $x$ and the initial population distribution $\mu_0$, 
thereby laying a solid foundation for further analysis 
of the value function's properties.
We employ the stochastic maximum principle and infinite-time forward-backward stochastic differential equations (FBSDEs)
to solve this game-theoretic problem and derive a representation of the equilibrium strategy via an elliptic master equation.

The theory of general nonlinear BSDEs was pioneered by Pardoux and Peng \cite{pardoux1990adapted,pardoux2005backward} in the
early 1990s, which is now a typical tool in stochastic optimization problems. 
We shall solve the mean field game problem  using the Pontryagin's maximum principle and infinite-time
FBSDEs.
The foundational work in \cite{peng2000infinite} proved the existence and uniqueness of solutions to
infinite-time FBSDEs, 
and subsequent work in \cite{shi2020forward} investigated a more general setting and established connections with quasilinear elliptic
partial differential equations (PDEs).
Recently,  \cite{ bayraktar2023solvability} extended this framework to  
McKean-Vlasov FBSDEs,
which play  crucial roles in obtaining  equilibrium solutions for mean field games.
In this paper, 
we partition all agents into a representative player and the
social equilibrium and characterize the equilibrium state  through the following system of infinite-time FBSDE:
\begin{equation}
  \label{eq: mv xi1}
  \begin{cases}
    \rmd X_{t}^{\xi} = b(X_{t}^{\xi}, \mathcal{L} _{X_{t}^{\xi}}, \hat{\alpha}(X_{t}^{\xi}, Y_{t}^{\xi})) \rmd t +\rmd B_{t}, \\
    \rmd Y_{t}^{\xi} = -\partial_x \mathcal{H} (X_{t}^{\xi},\mathcal{L} _{X_{t}^{\xi}}, \hat{\alpha}(X_{t}^{\xi}, Y_{t}^{\xi}), Y_{t}^{\xi}) \rmd t + Z_{t}^{\xi}\rmd B_{t}, \\
    X_0^{\xi}=\xi.
    \end{cases}
\end{equation}
\begin{equation}
  \label{eq: mv x1}
  \begin{cases}
    \rmd X_{t}^{x} = b(X_{t}^{x}, \mathcal{L} _{X_{t}^{\xi}}, \hat{\alpha}(X_{t}^{x}, Y_{t}^{x})) \rmd t +\rmd B_{t}, \\
    \rmd Y_{t}^{x} = -\partial_x \mathcal{H} (X_{t}^{x},\mathcal{L} _{X_{t}^{\xi}}, \hat{\alpha}(X_{t}^{x}, Y_{t}^{x}), Y_{t}^{x}) \rmd t + Z_{t}^{x}\rmd B_{t}, \\
    X_0^x=x.
    \end{cases}
\end{equation}
Here, $\mathcal{H}(x,\mu,a,y)\triangleq b(x,\mu,a)\cdot y+f(x,\mu,a)-rxy$ is the generalized Hamiltonian,
and $\hat{\alpha}(x,y) $ is the minimizer of $\mathcal{H}(x,\mu,a,y)$ with respect to $a$ when 
we assume that $\mathcal{H}$ is separable in variables $\mu$ and $a$.
It is worth emphasizing that we set the diffusion coefficient to 1 here, which makes both equations (\ref{eq: mv xi1}) 
and (\ref{eq: mv x1}) solvable. If the diffusion term in problem setup (\ref{eq: constraint}) is 
dependent on the state, distribution, 
and control, studying these two equations would be considerably more challenging.
Now, we consider the solutions to these two FBSDEs.
The solution $X_t^{\xi}$ to Equation (\ref{eq: mv xi1}) represents the population's state process, 
whose law corresponds to the  population distribution $\mu_t$.
In addition, the solution $X_t^x$ to Equation (\ref{eq: mv x1}) is the state process of the representative agent after solving the optimization problem. Notably, 
it exhibits the same feedback structure as the population's state process.
To further elucidate the relationship in Equation (\ref{eq: relationship}), 
we introduce  elliptic master equations.

First introduced by Lions in lectures \cite{lions2007cours}, the parabolic master equation appeared in the context
of the theory of mean field games. This is a time-dependent equation that bears profound connections with finite-time mean field game theory. 
Essentially, it describes a strategic interaction between a representative player and the collective environment. 
When the Nash equilibrium exists, the master equation provides a powerful tool to characterize the equilibrium cost and control pattern of this system. 
We refer the
reader to \cite{carmona2018probabilistic,cardaliaguet2019master,gangbo2022mean} for a comprehensive exposition on the subject.
In this paper, we propose  elliptic master equations, which explicitly characterize the feedback forms of both the representative player and the social equilibrium.
These equations take the following form:
\begin{equation}
\label{eq: master}
\begin{split}
  r U(x,\mu)=&H(x,\mu,\partial_xU(x,\mu))+\frac{1}{2}\partial_{xx}U(x,\mu)\\&+\tilde{\mathbb{E} }
  \left[\frac{1}{2}\partial_{\tilde{x}}\partial_\mu U(x,\mu,\tilde{\xi })+\partial_\mu U(x,\mu,\tilde{\xi })
  \partial_y H(\tilde{\xi },\mu,\partial_xU(\tilde{\xi},\mu))\right].
\end{split}
\end{equation}
Here, $\partial_x$ and $\partial_{xx}$ are standard spatial derivatives, 
$\partial_{\mu}$ and $\partial_{\tilde{x}\mu}$ are $W_2$-Wasserstein derivatives, $\tilde \xi$ is a random variable with law $\mu$ and $\tilde \dbE$ is the expectation with respect to its law. 
Under the assumption that the master equation (\ref{eq: master}) admits a solution with sufficient regularity, we derive the following representation for Equations (\ref{eq: mv xi1}) and (\ref{eq: mv x1}):
\begin{equation}
  Y_t^\xi=\partial_xU(X_{t}^{\xi}, \mathcal{L} _{X_{t}^{\xi}}),\quad
  Z_t^\xi=\partial_{xx}U(X_{t}^{\xi}, \mathcal{L} _{X_{t}^{\xi}}).
\end{equation}
\begin{equation}
  Y_t^x=\partial_xU(X_{t}^{x}, \mathcal{L} _{X_{t}^{\xi}}),\quad
  Z_t^x=\partial_{xx}U(X_{t}^{x}, \mathcal{L} _{X_{t}^{\xi}}).
\end{equation}
If $\tilde{b}(x,\mu)\triangleq b(x,\mu,\hat{\alpha}(x,\partial_xU(x,\mu)))$ is Lipschitz continuous in $(x,\mu)$,  
we have
\begin{equation}
  X_t^{x}|_{x=\xi}=X_t^{\xi},\quad t\in[0,T]
\end{equation}
for any fixed finite time $T$.
This is precisely the relationship expressed in Equation (\ref{eq: relationship}).

Alternatively, if we define $\mcv(x,\mu)\triangleq Y^{x,\xi}_0 $,
we demonstrate that $ \mcv$ serves as a viscosity solution to the distribution-dependent PDE below:
\begin{equation}
  \label{eq: intro-pa}
  \begin{split}
  r \mcu(x,\mu)=&\partial_xH(x,\mu,\mcu(x,\mu))+ \partial_yH(x,\mu,\mcu(x,\mu))\cdot \partial_x\mcu(x,\mu) +
  \frac{1}{2}\partial_{xx}\mcu(x,\mu)\\&+\tilde{\mathbb{E} }
  \left[\frac{1}{2}\partial_{\tilde{x}}\partial_\mu \mcu(x,\mu,\tilde{\xi })+\partial_\mu \mcu(x,\mu,\tilde{\xi })
  \partial_y H(\tilde{\xi },\mu,\mcu(\tilde{\xi},\mu))\right].
\end{split}
\end{equation}

This paper is organized as follows: in section \ref{sec:setting}, we present preliminaries for 
problems in this paper;
in section \ref{sec:assum}, we introduce infinite-time mean field games with discounting and the definition
of a Nash equilibrium; in section \ref{sec: fbsde} we characterize the equilibrium states through a system of infinite-time
FBSDEs;  in section \ref{sec:main results}, we introduce an elliptic master equation to provide a representation for the Nash equilibrium;
and in section \ref{sec: viscosity}, 
we prove the Yamada-Watanabe theorem  for infinite-time FBSDEs and 
provide a viscosity solution for distribution-dependent elliptic PDEs by virtue of the class of FBSDEs introduced in section \ref{sec: fbsde}.

\section{Preliminaries} 
\label{sec:setting}
\setcounter{equation}{0}

We will use a filtered probability space $(\Omega, \mathcal{F} ,\mathbb{P},\mathbb{F} )$
endowed
with a Brownian motion $B$.
Its filtration $\mathbb{F} \triangleq (\mathcal{F} _t)_{t\ge 0}$ is
augmented
by all $\mathbb{P}$-null sets and a sufficiently rich sub-$\sigma$-algebra $\mathcal{F}_0$ independent
of $B$ such that it
can support any measure on $ \mathbb{R} $ with finite second moment.

Let $( \O', \mathcal{F}' , \dbP',\dbF')$ be a copy of the filtered probability space $(\Omega, \mathcal{F} ,\mathbb{P},\mathbb{F})$
 with corresponding Brownian motion $B'$. We define the larger filtered probability space by
\begin{equation}
\tilde \O \triangleq \O\times \O' ,\q \tilde{ \mathcal{F}} \triangleq \mathcal{F}\otimes \mathcal{F}'\q \tilde\dbF = \{\tilde \cF_t\}_{t\ge 0} \triangleq \{\cF_t \otimes  \cF'_t\}_{t\ge 0},\q \tilde \dbP \triangleq \dbP\otimes \dbP',\q \tilde \dbE\triangleq \dbE^{\tilde \dbP}.
\end{equation}
Throughout the paper, we will use the probability space $(\Omega, \mathcal{F}, \mathbb{P}, \mathbb{F})$. However, when we deal with the distribution-dependent master equation, independent copies of random variables or processes are needed. Then, we will tacitly use their extensions to the larger space $(\tilde \O,\tilde{\mathcal{F}},  \tilde \dbP,\tilde \dbF)$.

Let $\cP \triangleq\cP(\dbR)$ be the set of all probability measures on $\mathbb R$ and let $\cP_p(p\ge1)$ denote the set of $\mu\in \cP$ with finite $p$-th moment.
For any sub-$\si$-field $\cG\subset \cF$ and $\mu\in \cP_p$, we define $\dbL^p(\cG)$ to be the set of $\dbR$-valued, $\cG$-measurable, and $p$-integrable random variables $\xi$ and $\dbL^p(\cG;\mu)$ to be the set of $\xi\in \dbL^p(\cG)$ such that the law $\cL_\xi=\mu$ .
For any $\mu,\nu\in \cP_p$, we define the $\cW_p$--Wasserstein distance between them as follows: 
\beaa
\cW_p(\mu, \nu) \triangleq \inf\Big\{\big(\dbE[|\xi-\eta|^q]\big)^{1/ q}: \mbox{for all $\xi\in \dbL^p(\cF; \mu)$, $\eta\in \dbL^p(\cF; \nu)$}\Big\}.
\eeaa

Due to our interest in discounted infinite-time mean field games, for any $K\in \mathbb{R} $, we denote by $L^2_K(t_0, \infty, \mathbb{R} )$
 the Hilbert space of all $\mathbb{R} -$valued adapted stochastic processes ($v_t$) starting from $t_0$ such that
\begin{equation}\mathbb{E}\left[\int_{t_0}^\infty e^{-Kt}|v_t|^2dt\right]<+\infty.\end{equation}
To simplify, we set $L^2_K \triangleq L^2_K(0, \infty,\mathbb{R} )$.
For each $\mathcal{F} _0$-measurable square-integrable random
variable $\xi$ , we consider the following infinite-time FBSDE:
\begin{equation}
  \label{eq: fbsde}
  \begin{cases}
    \rmd X_t=G(t,X_t,Y_t,\mathcal{L}_{X_t})dt+ \rmd B_t, \\
    \rmd Y_t=-F(t,X_t,Y_t,\mathcal{L}_{X_t})dt+Z_t\rmd B_t, \\
    X_0=\xi.
    \end{cases}
\end{equation}
Here, $G,F:\dbR^+\times \dbR^2\times \cP_2\to\dbR$ are two measurable functions that satisfy
 the following assumptions.
\begin{assum}
  \label{assum: fbsde}
  (i) There exists a positive constant $\ell$ such that for any $x,x',y,y'\in \mathbb{R}$ and $\mu,\mu'\in \cP_2$,
  \begin{equation}
    \begin{aligned}
      |G(t,x,y,\mu) & -G(t,x^{\prime},y^{\prime},\mu^{\prime})|+|F(t,x,y,\mu)-F(t,x^{\prime},y^{\prime},\mu^{\prime})| \\
       & \leq \ell(|x-x^{\prime}|+|y-y^{\prime}|+\mathcal{W}_{2}(\mu,\mu^{\prime})).\quad\mathrm{a.s.}
      \end{aligned}
  \end{equation}

(ii) There exist constants $0<K<2\kappa$ such that for any $t\geq0$ and any square-integrable 
random variables $X,X',Y,Y'$,
\begin{equation}
  \label{eq: condition}
  \begin{aligned}
    &\mathbb{E}\left[-K\hat{X}\hat{Y}-\hat{X}(F(t,U)-F(t,U^{\prime}))+\hat{Y}(B(t,U)-B(t,U^{\prime}))\right] \\
      &\leq-\kappa\mathbb{E}\left[\hat{X}^{2}+\hat{Y}^{2}\right],
    \end{aligned}
\end{equation}
where $ \hat{X}=X-X^{\prime},\hat{Y}=Y-Y^{\prime}$ and $U=(X,Y,\mathcal{L}_X),U'=(X',Y',\mathcal{L}_{X'}) .$
\end{assum}
The following lemma establishes the existence and uniqueness of a solution to the FBSDE(\ref{eq: fbsde}).
For a detailed proof, we refer the reader to \cite{bayraktar2023solvability} (Theorem 2.1).
\begin{lem}
  \label{lem: fbsde}
  Under Assumption \ref{assum: fbsde}, the FBSDE (\ref{eq: fbsde}) admits a unique solution in $L_K^2$.
\end{lem}

We introduce the Wasserstein space and the associated differential calculus. 
For a $\cW_2$-continuous function $U: \cP_2 \to \dbR$, its $\cW_2$-Wasserstein derivative\cite{carmona2018probabilistic}(also called the Lions-derivative) takes the form $\pa_\mu U: (\mu,\tilde x)\in \cP_2\times \dbR\to \dbR$ and satisfies
\bea
\label{pamu}
U(\cL_{\xi +  \eta}) - U(\mu) = \dbE\big[\langle \pa_\mu U(\mu, \xi), \eta \rangle \big] + o(\|\eta\|_2), \ \forall\ \xi\in\mathbb L^2(\mathcal{F};\mu),\eta\in\mathbb L^2(\mathcal{F}).
\eea
Let $C^0(\cP_2)$ denote the set of $\cW_2$-continuous functions $U:\cP_2\to\dbR$.
For $C^1(\cP_2)$, we define the space of functions $U\in C^0(\cP_2)$ such that $\pa_\mu U$ exists, is continuous on $ \cP_2\times \dbR$, and is uniquely determined by \reff{pamu}.
Let $C^{2,1}(\dbR\times\cP_2)$ denote the set of continuous functions $U:\dbR\times\cP_2\to\dbR$ such that $\pa_xU,\pa_{xx}U$ exist and are joint continuous on $\dbR\times \cP_2$ and $\pa_\mu U,\pa_{x\mu}U,\pa_{\tilde x\mu}U$ exist and are continuous on $\dbR\times\cP_2\times\dbR$.

Finally, we consider the space 
$\Th\triangleq [0, T]\times \dbR \times \cP_2$ for some $T>0$;
let $C^{1,2,1}(\Th)$ denote the set of  continuous functions $U: \Th\to \dbR$ that has the following continuous derivatives: 
$\pa_t U$, $\pa_x U$, $\pa_{xx} U$, $\pa_\mu U$, $  \pa_{x\mu} U$, $\pa_{\tilde x\mu} U.$ 
One crucial property of functions $U\in C^{1,2,1}(\Th)$ is that they satisfy the  It\^{o}'s  formula\cite{buckdahn2017mean,carmona2018probabilistic}. For $i=1,2$, let $\rmd X^i_t \triangleq b^i_t \rmd t + \si^i_t \rmd B_t ,$ where $b^i:[0,T]\times\Omega\to\mathbb R$ and $\sigma^i:[0,T]\times\Omega\to\mathbb R$ are $\dbF$-progressively 
measurable and bounded (for simplicity), and let $\rho_t\triangleq \cL_{X^2_t}$.
Suppose that for every compact subset $\cK\subset \dbR\times \cP_2$, it holds that
\begin{equation}
  \sup_{(t,x,\mu)\in [0,T]\times\cK} 
    \int_\dbR \left( \left\lvert \pa_\mu U(t,x,\mu,\tilde{x})\right\rvert^2 
    +\left\lvert \pa_{\tilde{x}\mu}U(t,x,\mu,\tilde{x})\right\rvert^2 \right) \rmd \mu(\tilde{x})
  <\infty.
\end{equation}
We have
\begin{equation}
  \label{eq:ito}
\begin{split}
  \rmd U(t, X^1_t, \rho_t) = & \Big[\pa_t U + \pa_x U\cdot b^1_t + \frac{1}{2} \pa_{xx} U (\si_t^1 )^2\Big](t, X^1_t, \rho_t) \rmd t \\
&+\Big(\tilde \dbE_{\cF_t}\big[\pa_\mu U(t,X^1_t,\rho_t,\tilde X^2_t) (\tilde b^{2}_t)+\frac{1}{2} \pa_{\tilde x}\pa_\mu U(t, X^1_t, \rho_t, \tilde X^2_t)(\tilde \si_t^2 )^2\big]\Big) \rmd t\\
&+\pa_xU(t,X^1_t,\rho_t)\si_t^1\rmd B_t  .
\end{split}
\end{equation}
Here,  $\tilde \dbE_{\cF_t}$ is the conditional expectations under $\tilde \dbP$
 given $\cF_t$ 
and the process $(\tilde X^2_t,\tilde b^{2}_t,\tilde\si_t^2)_{0\leq t\leq T}$
is a copy of the process $( X^2_t, b^{2}_t,\si_t^2)_{0\leq t\leq T}$, defined
on a copy $( \O', \mathcal{F}' , \dbP',\dbF')$ of the probability space $(\Omega, \mathcal{F} ,\mathbb{P},\mathbb{F} )$.

\section{Infinite-time mean field games with discounting} 
\label{sec:assum}
\setcounter{equation}{0}

In this section, we introduce  infinite-time mean field games with discounting.
Let $r>0$ represent the time discount factor and $A\subset \mathbb{R} $ be a convex control space.
Define $\mathcal{A} \triangleq L^2_r(0,\infty,A)$ to be  the space of all admissible controls,
and $b, f:\mathbb{R} \times \cP_2\times A\rightarrow\mathbb{R}  $ are two measurable functions.

We consider a population consisting of a continuum of players, where each individual player strategically 
interacts with the aggregate distribution formed by all other players to minimize their own cost.
Let $\mu_t$ denote the population distribution and $\xi\in \mathbb{L} ^2(\mathcal{F} _0)$ denote the initial
state with $\mathcal{L} _{\xi}=\mu_0$.
The state of the representative player with initial value $x$ is given by
\begin{equation}
  X_t^{x,\beta}=x+\int_{0}^{t}b(X_s^{x,\beta},\mu_t,\beta_s)\rmd s+B_t,
\end{equation}
where $\beta\in \mathcal{A} $ is the strategy that remains to be determined.
The representative player seeks to minimize the cost
\begin{equation}
  J^\mu(\beta)\triangleq\mathbb{E}\left[\int_0^\infty e^{-rt}f(X_t^{x,\beta},\mu_t,\beta_t)\rmd t\right].
\end{equation}
Assuming that $\beta^x\in \mathcal{A}$ minimizes the cost function, the state process of the representative player
is $X^{x,\beta^x}_t$.
Since the representative player can characterize the strategies of all players with the same initial state in the population, we assert that the following fundamental relationship must hold:
\begin{equation}
  \mu_t=\mathcal{L} _{X^{x,\beta^x}_t|_{x=\xi}}.
\end{equation}

To solve this mean field game problem, we work under the following assumptions.
\begin{assum}
  \label{assum: sde}
  (i) $b(x,\mu,a)$ is Lipschitz in $(x,\mu,a)$, and $f(x,\mu,a)$ is of at
  most quadratic growth in $(x,\mu,a)$. There exists a positive constant $\ell$ such that for
  any $\mu,\mu'\in \cP_2, x\in\mathbb{R} ,a\in A$, 
  \begin{equation}|b(x,\mu,a)-b(x,\mu^{\prime},a)|\leq \ell\mathcal{W}_2(\mu,\mu^{\prime}).\end{equation}

  (ii) There exists a constant $\lambda>\ell-r/2$ such that for any $a \in A, \mu \in \cP_2,
  x, x'\in \mathbb{R} $, it holds that
  \begin{equation}(x-x') \left( b(x,\mu,a)-b(x',\mu,a)\right)\leq -\lambda (x-x')^2.\end{equation}
\end{assum}
These assumptions jointly ensure that
(i) the state process remains confined within the $L^2_r$
  space and
(ii) the cost functional maintains integrability.

Then, we  partition all agents into two components: (i) the \emph{representative player}, who dynamically optimizes their strategy, and (ii) 
the \emph{social equilibrium} (or \emph{mean field}), characterizing the macroscopic state shared 
by the population.
The state of the social equilibrium is governed by the following stochastic differential equation (SDE):
\begin{equation}
  \label{eq: a state}
  X_t^{\xi,\alpha}=\xi+\int_0^t b( X_s^{\xi,\alpha}, \mathcal{L}_{ X_s^{\xi,\alpha}},\alpha_s)\mathrm{d}s+B_t,
\end{equation}
where $\xi\in \mathbb{L} ^2(\mathcal{F} _0)$ and $\alpha\in \mathcal{A} $. We note that by  assumption \ref{assum: sde}, 
this SDE has a unique strong solution in $L^2_r$; see \cite{bayraktar2023solvability} (Proposition 2.2) for more details.

The state of the representative player is governed by
\begin{equation}
  \label{eq: b state}
  X_t^{x,\beta}=x+\int_0^t b( X_s^{x,\beta}, \mathcal{L}_{ X_s^{\xi,\alpha}},\beta_s)\mathrm{d}s+B_t.
\end{equation}
Here, we also require their control $\beta\in \mathcal{A} .$

The representative player seeks to minimize the cost
\begin{equation}
  \begin{aligned}
    J\left(x,\xi;\alpha,\beta\right) =\mathbb{E}\bigg[\int_{0}^{+\infty}\rme^{-rt}f\big(X_{t}^{x,\beta},\mathcal{L}_{X_{t}^{\xi,\alpha}},\beta_t\big)\rmd t\bigg].
    \end{aligned}
\end{equation}
For any $(x,\xi)\in \dbR\times \mathbb{L} ^2(\mathcal{F} _0)$ and $\alpha\in \mathcal{A} $, we consider the infimum
\begin{equation}V(x,\xi;\alpha)\triangleq \inf_{\beta\in\mathcal{A}}J(x,\xi;\alpha,\beta).\end{equation}

\begin{defn}
  \label{def: ne}
  We say that a Lipschitz function ${\alpha}^*(x,\mu):\mathbb{R}\times \cP_2 \rightarrow \mathbb{R} $ 
  constitutes a discounted infinite-time  mean field Nash equilibrium for a given initial distribution $\mu_0$ if
  for any initial state $ \xi_0\in \mathbb{L} ^2(\mathcal{F} _0)$ with distribution $\mu_0$,
   the closed-loop controls $\alpha_s^{\xi_0}={\alpha}^*(X_s^{\xi_0,\alpha^{\xi_0}},\mathcal{L}_{X_s^{\xi_0,\alpha^{\xi_0}}} )$,$ \alpha_s^x={\alpha}^*(X_s^{x,\alpha^x},\mathcal{L}_{X_s^{\xi_0,\alpha^{\xi_0}}})$
  satisfy
  \begin{equation}
    \label{nash e}
    \alpha^x\in\arg\min_{\beta\in\mathcal{A} }J(x,\xi_0;\alpha^{\xi_0},\beta).
  \end{equation} 

  When the Nash equilibrium $\alpha^*$ exists, we have
\begin{equation}
  \xi,\xi^{\prime}\in\mathbb{L}^2(\mathcal{F}_0),\quad\mathcal{L}_{\xi^{\prime}}=
  \mathcal{L}_\xi\quad\Longrightarrow\quad  \mathcal{L}_{X_t^{\xi,\alpha^{\xi}}}   =\mathcal{L}_{X_t^{\xi',\alpha^{\xi'}}},\quad ~\mbox{for a.e. $t\ge 0$.}
\end{equation}
Therefore, we can define
\begin{equation}
  V(x,\mu)\triangleq J(x,\xi_0;\alpha^{\xi_0},\alpha^x), \quad\quad \xi_0\in \mathbb{L} ^2(\mathcal{F} _0,\mu).
\end{equation}
\end{defn}

In our framework, we separate a representative player from the population, who only needs to consider their  optimization 
problem starting from state $x$. 
This model decouples the micro-level agent  from the macro-level 
population distribution, enabling an interconnected analysis of their evolution. 
The equilibrium is characterized by two consistency conditions.
\begin{itemize}
  \item Individual Rationality: The representative player's optimal strategy is consistent with the perceived social equilibrium.
  \item Macro-consistency: The aggregate distribution generated by all players adopting this strategy must coincide with the posited social equilibrium.
\end{itemize}
The representative player's state evolution depends on both their state $X^{x,\beta}$ and the overall population distribution $\mathcal{L} _{X^{\xi,\alpha}}$.
Here, we would like to emphasize that since the number of players is large, any change by a representative player does not impact the measure flow $\mathcal{L} _{X^{\xi,\alpha}}$.
Under the existence assumption of the Nash equilibrium $\alpha^*$ specified in Definition \ref{def: ne}, 
the stochastic dynamics of both the population and  representative player are characterized by 
the following  SDE:
\begin{equation}
  \begin{cases}
  X_t^{\xi,\alpha^*}=\xi+\int_0^t b( X_s^{\xi,\alpha^*}, \mathcal{L}_{ X_s^{\xi,\alpha^*}},{\alpha}^*(X_s^{\xi,\alpha^*},\mathcal{L}_{X_s^{\xi,\alpha^*}} ))\mathrm{d}s+B_t, \\
X_t^{x,\alpha^*}=x+\int_0^t b( X_s^{x,\alpha^*}, \mathcal{L}_{ X_s^{\xi,\alpha^*}},{\alpha}^*(X_s^{x,\alpha^*},\mathcal{L}_{X_s^{\xi,\alpha^*}}))\mathrm{d}s+B_t.
\end{cases}
\end{equation}
 We further assume that $\tilde{b}(x,\mu)\triangleq b(x,\mu,\alpha^*(x,\mu))$ is Lipschitz continuous in $(x,\mu)$.
 For any fixed finite time $T$, see  \cite{carmona2015forward,buckdahn2017mean} , we have
\begin{equation}
  X_t^{x,\alpha^*}|_{x=\xi}=X_t^{\xi,\alpha^*},\quad t\in[0,T].
\end{equation}
This implies that every sample point from the initial population follows the same evolutionary dynamics 
as our hypothesized representative player, which justifies the mathematical validity of using a single representative 
player to characterize the behavior of all individuals in the population.
Moreover, when all agents in the population adopt the same strategy as the representative player, 
their collective behavior precisely generates the aggregate distribution $\mathcal{L} _{X_t^{\xi,\alpha^*}}$ derived from our solution. 
This justifies why we refer to $X_t^{\xi,\alpha^*}$ as the social equilibrium.

\section{Connection with infinite-time McKean-Vlasov FBSDEs}
\label{sec: fbsde}
\setcounter{equation}{0}

In this section, we employ the maximum principle to solve the optimization problem for the representative player
and then use  infinite-time McKean-Vlasov FBSDEs to construct the optimal strategy for the representative player
such that the controls of the representative player and the social equilibrium share the same feedback form.
Our derivation is based on the following key assumptions on $b,f$:

\begin{assum}
  \label{assum: b,f}
(i) $b(x,\mu,a)=b_0(x,\mu)+b_1(x,a)$ and $f(x,\mu,a)=f_0(x,\mu)+f_1(x,a)$ where $b_0,f_0$ are measurable functions on $\mathbb{R} \times\cP_2$ and
$b_1,f_1$ are measurable functions on $\mathbb{R} \times A$.

(ii) $b,f$ are differentiable with respect to $(x,a)$, and $\partial_ab, \partial_af$ are Lipschitz continuous in $(x,a)$.

(iii) $H_0(x,\mu,a,y)\triangleq b(x,\mu,a)\cdot y+f(x,\mu,a)$ is convex with respect to $(x,a)$. $\min_{a\in A} H_0(x,\mu,a,y)$
      has a unique minimizer $\hat{\alpha}(x,y)$ that is Lipschitz continuous in $(x,y)$.
\end{assum}

\subsection{Pontryagin's maximum principle}

Assuming that the state of the social equilibrium $X_t^{\xi}$ is given, we consider the optimization problem for the representative player,
whose state is given by 
\begin{equation}
  X_t^{x,\beta}=x+\int_0^t b( X_s^{x,\beta}, \mathcal{L}_{ X_s^{\xi}},\beta_s)\mathrm{d}s+B_t.
\end{equation}
The cost functional takes the form
\begin{equation}
  J(\beta)\triangleq\mathbb{E}\left[\int_0^\infty e^{-rt}f(X_t^{x,\beta},\mathcal{L}_{X_t^{\xi}},\beta_t)\rmd t\right],
\end{equation}
and the representative player wants to solve the minimization
problem
\begin{equation}\inf_{\beta\in\mathcal{A}}J(\beta).\end{equation}
Let us formally derive the maximum principle for the infinite-time control problem.
Suppose $\beta$ is an optimal control, choose another admissible control $\gamma$, and
denote by $X^{x,\beta+\epsilon\gamma}$
 the state trajectory corresponding to the control $\beta+\epsilon\gamma$
. Let
\begin{equation}R_t=\lim_{\epsilon\to0}\frac{X_t^{x,\beta+\epsilon\gamma}-X_t^{x,\beta}}{\epsilon}\end{equation}
be the variation process. Then, it can be shown that $R$ satisfies
\begin{equation}
  \begin{cases}
    \rmd R_t & =\left(\partial_xb( X_t^{x,\beta}, \mathcal{L}_{ X_t^{\xi}},\beta_t)\cdot R_t+\partial_ab(X_t^{x,\beta}, \mathcal{L}_{ X_t^{\xi}},\beta_t)\cdot\gamma_t\right)\rmd t, \\
    R_0 & =0.
    \end{cases}
\end{equation}

The function $\beta\to J(\beta)$ is G\^{a}teaux differentiable in the direction $\beta$ and its derivative
is given by
\begin{equation}
  \left.\frac{d}{d\epsilon}J(\beta+\epsilon\gamma)\right|_{\epsilon=0}=\mathbb{E}\left[\int_0^\infty e^{-rt}\left(\partial_xf( X_t^{x,\beta}, \mathcal{L}_{ X_t^{\xi}},\beta_t)\cdot R_t+\partial_af( X_t^{x,\beta}, \mathcal{L}_{ X_t^{\xi}},\beta_t)\cdot\gamma_t\right)dt\right].
\end{equation}
Define the generalized Hamiltonian
\begin{equation}
  \mathcal{H}(x,\mu,a,y)\triangleq b(x,\mu,a)\cdot y+f(x,\mu,a)-rxy,
\end{equation}
and introduce the adjoint process, which is determined by an infinite-time BSDE:
\begin{equation}
  \rmd Y_t^{x,\beta}=-\left(\partial_x\mathcal{H}(X_t^{x,\beta}, \mathcal{L}_{ X_t^{\xi}},\beta_t,Y_t^{x,\beta})\right)\rmd t+Z_t^{x,\beta}\rmd B_t.
\end{equation}

Applying It\^{o}'s  formula to the process $(e^{-rt}R_tY_t^{x,\beta})$ and by simple computation, we can deduce that
\begin{equation}
  \left.\frac{\rmd}{\rmd\epsilon}J(\beta+\epsilon\gamma)\right|_{\epsilon=0}=
  \mathbb{E}\left[\int_0^\infty e^{-rt}\partial_a\mathcal{H}(X_t^{x,\beta}, \mathcal{L}_{ X_t^{\xi}},\beta_t,Y_t^{x,\beta})\cdot\gamma_t\rmd t\right].
\end{equation}
Thus when $\beta$ is an optimal admissible control with the associated stochastic processes
$(X_t^{x,\beta} , Y_t^{x,\beta} , Z_t^{x,\beta} )$, it holds that
\begin{equation}
  \mathcal{H}(X_t^{x,\beta}, \mathcal{L}_{ X_t^{\xi}},\beta_t,Y_t^{x,\beta})=\min_{a\in A}\mathcal{H}(X_t^{x,\beta}, \mathcal{L}_{ X_t^{\xi}},a,Y_t^{x,\beta}).
\end{equation}
Recalling our convexity assumptions on $b,f$ in Assumption \ref{assum: b,f}, we know that the representative player's optimal control $\beta$ takes a feedback form, that is, 
\begin{equation}
  \beta_t=\hat{\alpha}(X_t^{x,\beta},Y_t^{x,\beta}).
\end{equation}

Now, we consider the following two McKean-Vlasov FBSDEs:
\begin{equation}
  \label{eq: mv xi}
  \begin{cases}
    \rmd X_{t}^{\xi} =  b(X_{t}^{\xi}, \mathcal{L} _{X_{t}^{\xi}}, \hat{\alpha}(X_{t}^{\xi}, Y_{t}^{\xi})) \rmd t +\rmd B_{t}, \\
    \rmd Y_{t}^{\xi} = -\partial_x \mathcal{H} (X_{t}^{\xi},\mathcal{L} _{X_{t}^{\xi}}, \hat{\alpha}(X_{t}^{\xi}, Y_{t}^{\xi}), Y_{t}^{\xi}) \rmd t + Z_{t}^{\xi}\rmd B_{t}, \\
    X_0^{\xi}=\xi.
    \end{cases}
\end{equation}
\begin{equation}
  \label{eq: mv x}
  \begin{cases}
    \rmd X_{t}^{x} = b(X_{t}^{x}, \mathcal{L} _{X_{t}^{\xi}}, \hat{\alpha}(X_{t}^{x}, Y_{t}^{x})) \rmd t +\rmd B_{t}, \\
    \rmd Y_{t}^{x} = -\partial_x \mathcal{H} (X_{t}^{x},\mathcal{L} _{X_{t}^{\xi}}, \hat{\alpha}(X_{t}^{x}, Y_{t}^{x}), Y_{t}^{x}) \rmd t + Z_{t}^{x}\rmd B_{t}, \\
    X_0^x=x.
    \end{cases}
\end{equation}
Here, (\ref{eq: mv xi}) and (\ref{eq: mv x}) denote the state processes of the social equilibrium and the representative player, respectively. 
Observe that their admissible controls both take the identical feedback form $\hat{\alpha}(x,y)$. 
We shall prove that when the social equilibrium employs this feedback control, the representative player's loss function 
is minimized if they use the same feedback form, thereby constituting a Nash equilibrium. 

\begin{thm}
  Let $(b, f )$ be differentiable in $(x, a)$ and $\mathcal{H} $ be
convex in $(x, a)$. Suppose that $\hat{\alpha} $ is Lipschitz continuous and that both (\ref{eq: mv xi})  and (\ref{eq: mv x}) admit unique strong solutions in $L^2_r$.
If we denote $\hat{\alpha}(X_{t}^{x}, Y_{t}^{x})$ as $\alpha^*_t$, which is an admissible control in $\mathcal{A}$, then we have
\begin{equation}
  J({\alpha}^*)=\min_{\beta\in \mathcal{A}}J(\beta).
\end{equation}
\end{thm}
\proof
For an arbitrary
admissible control $\beta\in \mathcal{A} $ and its associated process 
\begin{equation}
  X_t^{x,\beta}=x+\int_0^t b( X_s^{x,\beta}, \mathcal{L}_{ X_s^{\xi}},\beta_s)\mathrm{d}s+B_t,
\end{equation}
we have 
\begin{equation}
  \label{eq: 4.1}
  \begin{aligned}
    J(\alpha^*)-J(\beta)= & \mathbb{E}\left[\int_{0}^{\infty}e^{-rt}\left(\mathcal{H}(X_{t}^{x},\mathcal{L} _{X_{t}^{\xi}}, \alpha^*_t, Y_{t}^{x})-\mathcal{H}(X_{t}^{x,\beta},\mathcal{L} _{X_{t}^{\xi}}, \beta_t, Y_{t}^{x})\right)\rmd t\right] \\
     & -\mathbb{E}\left[\int_0^\infty e^{-rt}\left(b(X_{t}^{x},\mathcal{L} _{X_{t}^{\xi}}, \alpha^*_t)-b(X_{t}^{x,\beta},\mathcal{L} _{X_{t}^{\xi}}, \beta_t)\right)\cdot Y_t^x\rmd t\right] \\
     & +r\mathbb{E}\left[\int_0^\infty e^{-rt}(X_t^x-X_t^{x,\beta})\cdot Y_t^x\rmd t\right]. 
    \end{aligned}
\end{equation}
Since $X_t^x,X_t^{x,\beta}, Y_t^x$ all belong to $L^2_r$, we can find a  sequence of $T_i\rightarrow \infty$
such that 
\begin{equation}
  \mathbb{E}\left[e^{-rT_{i}}(X_{T_{i}}^x-X_{T_{i}}^{x,\beta})\cdot Y_{T_{i}}^x\right]\to0.
\end{equation}
Applying It\^{o}'s formula to $e^{-rT_{i}}(X_{T_{i}}^x-X_{T_{i}}^{x,\beta})\cdot Y_{T_{i}}^x$
and letting $T_i\to \infty$, we obtain that
\begin{equation}
  \label{eq: 4.2}
  \begin{aligned}
    & \mathbb{E}\left[\int_{0}^{\infty}e^{-rt}(X_{t}^x-X_{t}^{x,\beta})\left(\partial_{x}\mathcal{H}(X_{t}^{x},\mathcal{L} _{X_{t}^{\xi}}, \alpha^*_t, Y_{t}^{x})\right)\rmd t\right] \\
     =&\mathbb{E}\left[\int_{0}^{\infty}e^{-rt}\left(-r(X_{t}^x-X_{t}^{x,\beta})+b(X_{t}^{x},\mathcal{L} _{X_{t}^{\xi}}, \alpha^*_t)-b(X_{t}^{x,\beta},\mathcal{L} _{X_{t}^{\xi}}, \beta_t)\right)\cdot Y_{t}^x\rmd t\right].
   \end{aligned}
\end{equation}
According to the convexity and differentiability of $\mathcal{H} $, we have
\begin{equation}
  \label{eq: 4.3}
  \begin{aligned}
    & \mathcal{H} (X_{t}^{x,\beta},\mathcal{L} _{X_{t}^{\xi}}, \beta_t, Y_{t}^{x})-\mathcal{H}(X_{t}^{x},\mathcal{L} _{X_{t}^{\xi}}, \alpha^*_t, Y_{t}^{x}) \\
     \geq&(X_{t}^{x,\beta}-X_{t}^x)\cdot\partial_{x}\mathcal{H}(X_{t}^{x},\mathcal{L} _{X_{t}^{\xi}}, \alpha^*_t, Y_{t}^{x})
     +(\beta_t-{\alpha}^*_t)\cdot\partial_a\mathcal{H}(X_{t}^{x},\mathcal{L} _{X_{t}^{\xi}}, \alpha^*_t, Y_{t}^{x}) \\
     =& (X_{t}^{x,\beta}-X_{t}^x)\cdot\partial_{x}\mathcal{H}(X_{t}^{x},\mathcal{L} _{X_{t}^{\xi}}, \alpha^*_t, Y_{t}^{x}).
   \end{aligned}
\end{equation}
The last equality follows from the fact that $\alpha^*_t=\hat{a}(X_{t}^{x}, Y_{t}^{x})$ is the minimizer of 
$\min_{a\in A}\mathcal{H}(X_{t}^{x},\mathcal{L} _{X_{t}^{\xi}}, a, Y_{t}^{x})$.
Combining equations (\ref{eq: 4.1}), (\ref{eq: 4.2}), and (\ref{eq: 4.3}), we obtain
\begin{equation}
  J(\alpha^*)-J(\beta)\leq 0
\end{equation}
for all admissible controls $\beta$. Thus, we complete the proof.
\qed

\subsection{Solvability of mean field game FBSDEs}
In this subsection, we will find sufficient conditions for the existence
and uniqueness of solutions to (\ref{eq: mv xi}) and (\ref{eq: mv x}).
Considering the linear case, we assume that $b(x,\mu,a)=b_1x+b_2\bar{\mu}+b_3a$ and $f(x,\mu,a)=f_0(x,\mu)+f_1(x,a)$,
where $\bar{\mu}$ is the mean value of the probability measure $\mu$ and $b_1,b_2,b_3,b_4$ are constants.

We require a technical lemma about the Lipschitz and convex properties of the minimizer $\hat{\alpha}$, the detailed proof of which can be found in \cite{bayraktar2023solvability} (Lemma 3.1).
\begin{lem}
  \label{lemma: alpha}
  Suppose $f_1$ is once continuously
  differentiable in $(x, a)$, $\partial_af_1$ is $\ell$-Lipschitz in  $x$ and $f_1$ is $\eta$-convex in $a$, which means that
  \begin{equation}
    f_1(x,a')-f_1(x,a)-(a'-a)\cdot \partial_af_1(x,a)\geq \eta|a'-a|^2,\quad ~\mbox{for all $x\in\mathbb{R} $.}
  \end{equation}
  Then, it holds that 
  \begin{equation}
    \label{eq: lip alpha}
    |\hat{\alpha}(x,y)-\hat{\alpha}(x^{\prime},y^{\prime})|\leq\frac{\ell}{2\eta}|x^{\prime}-x|+\frac{|b_3|}{2\eta}|y^{\prime}-y|,
  \end{equation}
  and for some $a_0\in A$,
  \begin{equation}
    |\hat{\alpha}(x,y)|\leq\eta^{-1}(|\partial_af_1(x,a_0)|+|b_2y|)+|a_0|.
  \end{equation}
  Furthermore, if $A=\mathbb{R} $ and $\partial_a f$ is  $\zeta$-Lipschitz in $a$, it follows that
  \begin{equation}
    \label{eq: lip yalpha}
    b_3(y^{\prime}-y)\cdot\left(\hat{\alpha}(x,y^{\prime})-\hat{\alpha}(x,y)\right)\leq-\frac{2b_3^2\eta}{\zeta^2}(y^{\prime}-y)^2.
  \end{equation}
\end{lem}

\begin{thm}
  \label{thm: solvable}
  Let $b(x,\mu,a)=b_1x+b_2\bar{\mu}+b_3a$ and $f(x,\mu,a)=f_0(x,\mu)+f_1(x,a)$. Under the following conditions,
Assumptions \ref{assum: sde} and \ref{assum: b,f} are satisfied and both (\ref{eq: mv xi})  and (\ref{eq: mv x}) admit unique strong solutions in $L^2_r$.

(i) There exists a positive constant $k$ such
that $|b_2| \leq k$ and $- b_1 \geq k - \frac{r}{2}$. $f_0(x,\mu)$ is once continuously differentiable in $x$
and of at most quadratic growth in $(x,\mu)$.
 $f_1(x,a)$ is once continuously differentiable
and of at most quadratic growth in $(x,a)$.

(ii) $f_0(x,\mu)$ is convex in $x$, and  there exists a positive constant $b_4>0$ such that
\begin{equation}
  |\partial_x f_0(x,\mu)-\partial_x f_0(x',\mu')|<b_4\left( |x-x'|+ \cW_2(\mu,\mu')  \right).
\end{equation}

(iii) There exist some positive constants $\eta,\iota$ such that the following convexity condition
holds:
\begin{equation}
  \begin{aligned}
    & f_1(x^{\prime},a^{\prime})-f_1(x,a)-\partial_{(x,a)}f_1(x,a)\cdot(x^{\prime}-x,a^{\prime}-a) \\
     \geq&\iota(x^{\prime}-x)^2+\eta(a^{\prime}-a)^2.
   \end{aligned}
\end{equation}

(iv) There exist some positive constants $\zeta, \ell$ such that $\partial_af_1$ is $\ell$-Lipschitz in  $x$ and $\zeta$-Lipschitz
in $a$. $\partial_xf_1$ is Lipschitz continuous in $(x,a)$.

(v) $A=\mathbb{R} $, and it holds that
\begin{equation}
  \label{eq: require}
  \min\left\{2\iota-\frac{\ell^2}{2\eta}-\frac{b_3\ell}{2\eta}-\frac{|b_2|}{2}-2b_4,
  \frac{2b_{3}^{2}\eta}{\zeta^{2}}-\frac{b_3\ell}{2\eta}-\frac{|b_2|}{2}\right\}>\frac{r}{2}.
\end{equation}

\end{thm}

\begin{rem}
  \label{rm: require detaile}
  If we set $f_0(x,\mu)=b_4x\bar{\mu}, f_1(x,a)=Ax^2+Ca^2, A>0, C>0$, we have $\iota=A,\eta=C,\zeta=2C,\ell=0$.
  Then, the requirement in (\ref{eq: require}) becomes
  \begin{equation}
    \label{eq: require detaile}
    \min\left\{2A-\frac{|b_2|}{2}-|b_4|,
  \frac{b_3^2}{2C}-\frac{|b_2|}{2}\right\}>\frac{r}{2}.
  \end{equation}
  Fixing $C,b_2,b_4,$ we take a large $b_3$ such that $\frac{b_3^2}{2C}-\frac{|b_2|}{2}$ is greater than $r/2$.
  Then we choose a sufficiently large $A$ such that $2A-\frac{|b_2|}{2}-|b_4|$ exceeds $r/2$.
  This construction satisfies all the required conditions.
\end{rem}

\proof
It is clear that Assumption \ref{assum: sde} is satisfied. And by Lemma \ref{lemma: alpha},
$\hat{\alpha}$ is Lipschitz. In addition, $\mathcal{H}(x,\mu,a,y)=(b_1x+b_2\bar{\mu}+b_3a)\cdot y + f_0(x,\mu) +f_1(x,a) $ is convex in $(x,a)$ according to the assumption on $f_0, f_1$;
thus, Assumption \ref{assum: b,f} is satisfied.
To prove that the FBSDEs (\ref{eq: mv xi})  and (\ref{eq: mv x}) admit unique strong solutions in $L^2_r$,
the only condition that remains to be verified is (\ref{eq: condition}).

We consider the FBSDE (\ref{eq: mv xi}) and set 
\begin{equation}
  \begin{aligned}
    & B(t,x,y,\mu)=b_1x+b_2\overline{\mu}+b_{3}\hat{\alpha}(x,y), \\
    & F(t,x,y,\mu)=b_{1}y+ \partial_x f_0(x,\mu) +\partial_{x}f_1(x,\hat{\alpha}(x,y))-ry .
   \end{aligned}
\end{equation}
Take four arbitrary square-integrable random variables $X,Y,X',Y'$.
Define $\hat{X}=X-X^{\prime},\hat{Y}=Y-Y^{\prime}$ and $U=(X,Y,\mathcal{L}_{X}),U^{\prime}=(X^{\prime},Y^{\prime},\mathcal{L}_{X'})$ .
We have
\begin{equation}
  \begin{aligned}
    &-r\hat{X}\hat{Y} -\hat{X}[F(t,U)-F(t,U^{\prime})]+\hat{Y}[B(t,U)-B(t,U^{\prime})] \\
   =&-r\hat{X}\hat{Y}-\hat{X}\left( (b_1-r)\hat{Y}+ \partial_x f_0(X,\mathcal{L}_{X})- \partial_x f_0(X',\mathcal{L}_{X'})   +\partial_xf_1(X,\hat{\alpha}(X,Y))- \partial_xf_1(X',\hat{\alpha}(X',Y')) \right) \\
     & + \hat{Y}\left(b_1\hat{X}+ {b}_{2}\mathbb{E}[\hat{X}]+b_{3}(\hat{\alpha}(X,Y)-\hat{\alpha}(X^{\prime},Y^{\prime}))\right)\\
  \le&-\hat{X}\left( \partial_xf_1(X,\hat{\alpha}(X,Y))- \partial_xf_1(X',\hat{\alpha}(X',Y')) \right)+b_3 \hat{Y}\left(\hat{\alpha}(X,Y)-\hat{\alpha}(X^{\prime},Y^{\prime})\right)\\
   &+b_2\hat{Y}\mathbb{E}[\hat{X}]+b_4|\hat{X}|\cdot \left(|\hat{X}| +  (\mathbb{E}[\hat{X}^2])^{1/2}\right).
    \end{aligned}
\end{equation}

Since $f_1$ is $\iota$-convex in $x$, we have
\begin{equation}
  \left[ \partial_xf_1(x',a)- \partial_xf_1(x,a) \right](x'-x)\geq 2\iota(x'-x)^2.
\end{equation}
Moreover, $\partial_xf_1$ is $l$-Lipschitz in $a$ and $\hat{\alpha}$ satisfies (\ref{eq: lip alpha}), we have that
\begin{equation}
  \begin{aligned}
   & -\hat{X}\left( \partial_xf_1(X,\hat{\alpha}(X,Y))- \partial_xf_1(X',\hat{\alpha}(X,Y)) \right) \\
   =&-\hat{X}\left( \partial_xf_1(X,\hat{\alpha}(X,Y))- \partial_xf_1(X',\hat{\alpha}(X,Y)) \right) \\
     & -\hat{X}\left( \partial_xf_1(X',\hat{\alpha}(X,Y))- \partial_xf_1(X',\hat{\alpha}(X',Y')) \right)\\
  \leq&-2\iota \hat{X}^2+l|\hat{X}|(\frac{l}{2\eta}| \hat{X}|+\frac{|b_3|}{2\eta}|\hat{Y}|).
    \end{aligned}
\end{equation}
From Lemma \ref{lemma: alpha}, $\hat{\alpha}$ satisfies (\ref{eq: lip yalpha}), it follows that
\begin{equation}
  \begin{aligned}
    &b_3 \hat{Y}\left(\hat{\alpha}(X,Y)-\hat{\alpha}(X^{\prime},Y^{\prime})\right) \\
 =&  b_3 \hat{Y}\left(\hat{\alpha}(X,Y)-\hat{\alpha}(X,Y^{\prime})\right)+b_3 \hat{Y}\left(\hat{\alpha}(X,Y')-\hat{\alpha}(X^{\prime},Y^{\prime})\right) \\
     & \leq-\frac{2b_{3}^{2}\eta}{\zeta^{2}}\hat{Y}^{2}+|b_3\hat{Y}|\left(\frac{l}{2\eta}|\hat{X}|\right).
    \end{aligned}
\end{equation}
By applying elementary estimates, we derive
\begin{equation}
  \begin{aligned}
    &\mathbb{E} \left[-r\hat{X}\hat{Y} -\hat{X}(F(t,U)-F(t,U^{\prime}))+\hat{Y}(B(t,U)-B(t,U^{\prime})) \right]\\
   \leq& (-2\iota+\frac{l^2}{2\eta}+\frac{b_3l}{2\eta}+\frac{|b_2|}{2}+2b_4)\mathbb{E}[\hat{X}^2]
   +(-\frac{2b_{3}^{2}\eta}{\zeta^{2}}+\frac{b_3l}{2\eta}+\frac{|b_2|}{2})\mathbb{E}[\hat{Y}^2]\\
   < & -\frac{r}{2}\mathbb{E}[\hat{X}^2 +\hat{Y}^2  ].
    \end{aligned}
\end{equation}

Now we know (\ref{eq: mv xi}) admits a unique solution $X_t^{\xi}\in L^2_r$ and set $\bar{\mu}_t=\mathbb{E}[X_t^{\xi} ]$
For the FBSDE(\ref{eq: mv x}), we set
\begin{equation}
  \begin{aligned}
    & B'(t,x,y)=b_1x+b_2\overline{\mu}_t+b_{3}\hat{\alpha}(x,y), \\
    & F'(t,x,y)=b_{1}y+\partial_x f_0(x,\mu_t)+\partial_{x}f_1(x,\hat{\alpha}(x,y))-ry .
   \end{aligned}
\end{equation}
Following the identical analytical procedure, we have 
\begin{equation}
  \begin{aligned}
    &\mathbb{E} \left[-r\hat{X}\hat{Y} -\hat{X}(F'(t,X,Y)-F'(t,X',Y'))+\hat{Y}(B'(t,X,Y)-B'(t,X',Y')) \right]\\
   \leq& (-2\iota+\frac{l^2}{2\eta}+\frac{b_3l}{2\eta}+b_4)\mathbb{E}[\hat{X}^2]
   +(-\frac{2b_{3}^{2}\eta}{\zeta^{2}}+\frac{b_3l}{2\eta})\mathbb{E}[\hat{Y}^2]\\
   < & -\frac{r}{2}\mathbb{E}[\hat{X}^2 +\hat{Y}^2  ].
    \end{aligned}
\end{equation}
By Lemma \ref{lem: fbsde}, we know both FBSDEs (\ref{eq: mv xi})  and (\ref{eq: mv x}) admit unique strong solutions in $L^2_r$.
\qed

\section{Master equation representation} 
\label{sec:main results}
\setcounter{equation}{0}

While we have derived a Nash equilibrium solution through  FBSDEs that yields identical feedback forms 
for both the representative player and social equilibrium, this feedback structure differs from our previously 
defined formulation in Definition \ref{def: ne}. In this section, we shall establish an alternative representation of the Nash equilibrium 
using classical solutions to the elliptic master equation (\ref{eq: master}).

Following Assumption \ref{assum: b,f}, we define 
\begin{equation}
  \label{eq: hamilton}
  H(x,\mu,y)=H_0(x,\mu,\hat{\alpha}(x,y),y).
\end{equation}
Through the assumptions on $b,f$, we can easily deduce the relationship
\begin{equation}
  \partial_yH(x,\mu,y)=b(x,\mu,\hat{\alpha}(x,y)).
\end{equation}

Assume that the master equation ($\ref{eq: master}$) has a classical solution $U(x,\mu)\in C^{2,1}(\dbR\times\cP_2)$ with $H$ defined above
such that $\pa_xU(x,\mu)\in C^{2,1}(\dbR\times\cP_2)$ satisfying PDE (\ref{eq: intro-pa}) and $\hat{b}(x,\mu)=b(x,\mu,\hat{\alpha}(x,\partial_xU(x,\mu)))$ is Lipschitz continuous in $(x,\mu)$.
For any $\xi\in \mathbb{L} ^2(\mathcal{F} _0)$, suppose the following SDE admits a unique solution in $L^2_r$:
\begin{equation}
  \label{eq: master xi state}
    \mathcal{X}_{t}^{\xi} = \xi + \int_0^t b\left( \mathcal{X}_{s}^{\xi}, \mathcal{L}_{\mathcal{X}_{s}^{\xi}}, \hat{\alpha}\left( \mathcal{X}_{s}^{\xi}, \partial_x U\left( \mathcal{X}_{s}^{\xi}, \mathcal{L}_{\mathcal{X}_{s}^{\xi}} \right) \right) \right) \mathrm{d}s + B_{t}.
\end{equation}
Then, we take $\mathcal{Y} _t^\xi=\partial_xU(\mathcal{X} _{t}^{\xi}, \mathcal{L} _{\mathcal{X} _{t}^{\xi}})$. 
By differentiating both sides of the master equation (\ref{eq: master}) with respect to $x$ and applying It\^{o}'s formula, we obtain
\begin{equation}
  \begin{split}
\rmd \mathcal{Y} _t^\xi =& \Bigg[ 
    \partial_{xx}U\left( \mathcal{X}_{t}^{\xi}, \mathcal{L}_{\mathcal{X}_{t}^{\xi}} \right) \cdot b\left( \mathcal{X}_{t}^{\xi}, \mathcal{L}_{\mathcal{X}_{t}^{\xi}}, \hat{\alpha}\left( \mathcal{X}_{t}^{\xi}, \partial_xU\left( \mathcal{X}_{t}^{\xi}, \mathcal{L}_{\mathcal{X}_{t}^{\xi}} \right) \right) \right) \\ 
    &+ \frac{1}{2}\partial_{xx}\partial_xU\left( \mathcal{X}_{t}^{\xi}, \mathcal{L}_{\mathcal{X}_{t}^{\xi}} \right) 
    + \tilde{\mathbb{E}}_{\mathcal{F}_t} \Bigg[ 
        \frac{1}{2}\partial_{\tilde{x}}\partial_{\mu}\partial_xU\left( \mathcal{X}_{t}^{\xi}, \mathcal{L}_{\mathcal{X}_{t}^{\xi}}, \tilde{\mathcal{X}}^{\xi}_t \right) \\
    &+ \partial_yH\left( \tilde{\mathcal{X}}_t^{\xi}, \mathcal{L}_{\mathcal{X}_{t}^{\xi}}, \partial_xU\left( \tilde{\mathcal{X}}_t^{\xi}, \mathcal{L}_{\mathcal{X}_{t}^{\xi}} \right) \right) \cdot
        \partial_{\mu}\partial_xU\left( \mathcal{X}^{\xi}_t, \mathcal{L}_{\mathcal{X}_{t}^{\xi}}, \tilde{\mathcal{X}}^{\xi}_t \right) 
    \Bigg] \Bigg] \rmd t \\
    &+ \partial_{xx}U\left( \mathcal{X}_{t}^{\xi}, \mathcal{L}_{\mathcal{X}_{t}^{\xi}} \right)  \rmd B_t \\
    =& \left( r \partial_{x}U\left( \mathcal{X}_{t}^{\xi}, \mathcal{L}_{\mathcal{X}_{t}^{\xi}} \right) - \partial_x H\left( \mathcal{X}_{t}^{\xi}, \mathcal{L}_{\mathcal{X}_{t}^{\xi}}, \partial_{x}U\left( \mathcal{X}_{t}^{\xi}, \mathcal{L}_{\mathcal{X}_{t}^{\xi}} \right) \right) \right) \rmd t 
    + \partial_{xx}U\left( \mathcal{X}_{t}^{\xi}, \mathcal{L}_{\mathcal{X}_{t}^{\xi}} \right)  \rmd B_t.
  \end{split}
\end{equation}
By comparing it with (\ref{eq: mv xi}), we derive the following relationship for social equilibrium:
\begin{equation}
  Y_t^\xi=\partial_xU(X_{t}^{\xi}, \mathcal{L} _{X_{t}^{\xi}}),\quad
  Z_t^\xi=\partial_{xx}U(X_{t}^{\xi}, \mathcal{L} _{X_{t}^{\xi}}).
\end{equation}
Applying the same argument to  (\ref{eq: mv x}), we obtain:
\begin{equation}
  Y_t^x=\partial_xU(X_{t}^{x}, \mathcal{L} _{X_{t}^{\xi}}),\quad
  Z_t^x=\partial_{xx}U(X_{t}^{x}, \mathcal{L} _{X_{t}^{\xi}}).
\end{equation}
This demonstrates that both the representative player and social equilibrium employ the same closed-loop control
\begin{equation}
  \label{eq: ne}
  \alpha^*(x,\mu)=\hat{\alpha}(x,\partial_xU(x,\mu)).
\end{equation}

We now revisit the mean field games through the master equation. Let $\xi\in \mathbb{L} ^2(\mathcal{F} _0)$ 
be the initial state with distribution $\mu$. We prove that under the assumption that the master equation admits 
a classical solution, 
the feedback control defined by (\ref{eq: ne}) constitutes a Nash equilibrium. Moreover, the solution to the master equation 
is precisely the value function of the representative player.

\begin{thm}
  Assume that the master equation (\ref{eq: master}) admits a classical solution $U(x,\mu)\in C^{2,1}(\dbR\times\cP_2)$, which is of at most quadratic
growth,
  and the following SDEs admit unique solutions in $L^2_r$ for each $x\in\mathbb{R} $ and $\xi\in\mathbb{L} ^2(\mathcal{F} _0;\mu)$
  \begin{equation}
  \begin{cases}
     X_{t}^{\xi} = \xi +\int_0^t b(X_{s}^{\xi}, \mathcal{L} _{X_{s}^{\xi}}, \hat{\alpha}(X_{s}^{\xi}, \partial_xU(X_{s}^{\xi}, \mathcal{L} _{X_{s}^{\xi}}))) \rmd t + B_{t},\\
X_{t}^{x} =x +\int_0^t b(X_{s}^{x}, \mathcal{L} _{X_{s}^{x}}, \hat{\alpha}(X_{s}^{x}, \partial_xU(X_{s}^{x}, \mathcal{L} _{X_{s}^{x}}))) \rmd t + B_{t}.
    \end{cases}
\end{equation}
  Then, for every initial distribution $\mu_0$, $\alpha^*(x,\mu)=\hat{\alpha}(x,\partial_xU(x,\mu))$
  is a Nash equilibrium, and $U(x,\mu_0)$ is exactly the value function given $\alpha^*$. 
\end{thm}
\proof 
Given $\alpha^*$ and initial $ \xi\in\mathbb{L} ^2(\mathcal{F} _0;\mu_0)$, the state of social equilibrium is governed by the following SDE:
\begin{equation}
  X_t^{\xi}=\xi+\int_0^t\partial_yH(X_s^{\xi},\rho_s,\partial_xU(X_s^{\xi},\rho_s))\rmd s+B_t,\quad\quad \rho_s\triangleq \mathcal{L} _{X_s^{\xi}}.
\end{equation}
We note that this SDE is the same as (\ref{eq: master xi state}).
Meanwhile, the state of the representative player is governed by  
\begin{equation}
  X_t^{x}=x+\int_0^t b(X_s^x,\rho_s,\beta_s)\mathrm{d}s+B_t,
\end{equation}
where $\beta\in \mathcal{A}$ remains to be determined.
Applying It\^{o}'s formula to $\rme^{-r t}U(X^x_t,\rho_t)$ yields
\begin{equation}
  \begin{split}
    \mathbb{E} \rme^{-r T}U(X^x_T,\rho_T)=&U(x,\mu_0)+\mathbb{E} \int_0^T\big[
      -r\rme^{-r t}U(X^x_t,\rho_t)
      +\rme^{-r t} \partial_x U(X^x_t,\rho_t)\cdot b(X_t^x,\rho_t,\beta_t)\\&
      +\frac{1}{2}\rme^{-r t}\partial_{xx}U(X^x_t,\rho_t)
      +\rme^{-r t}\tilde{\mathbb{E} }_{\mathcal{F}_t}\big[
        \frac{1}{2}\partial_{\tilde{x}}\partial_{\mu}U(X^x_t,\rho_t,\tilde{X}^{\xi}_t)
        \\&+\partial_yH(\tilde{X}_t^{\xi},\rho_t,\partial_xU(\tilde{X}_t^{\xi},\rho_t))\cdot
        \partial_{\mu}U(X^x_t,\rho_t,\tilde{X}^{\xi}_t)
      \big] \big]\rmd t\\
      =&U(x,\mu_0)+\mathbb{E} \int_0^T \rme^{-r t}\left[
        \partial_xU(X^x_t,\rho_t)\cdot b(X_t^x,\rho_t,\beta_t)-H(X^x_t,\rho_t,\partial_xU(X^x_t,\rho_t))
      \right]\rmd t\\
      \ge&U(x,\mu_0)-\mathbb{E} \int_0^T \rme^{-r t}f(X^x_t,\rho_t,\beta_t)\rmd t .   
  \end{split}
\end{equation}
Since $U$ is of at most quadratic
growth, taking the limit as $T\rightarrow +\infty$ yields
\begin{equation}
  U(x,\mu_0)\le\mathbb{E} \int_0^{+\infty} \rme^{-r t}f(X^x_t,\rho_t,\beta_t)\rmd t
\end{equation}
for every feasible control $\beta$, and the equality holds when $\beta_t=\hat{\alpha}(X_t^{x},\partial_xU(X_t^{x}, \rho_t))$.
This shows that $\alpha^*(x,\mu)=\hat{\alpha}(x,\partial_xU(x,\mu))$
  is a Nash equilibrium. In addition, we have
  \begin{equation}
    U(x,\mu_0)=\mathbb{E} \int_0^{+\infty} \rme^{-r t}f(X^x_t,\rho_t,\hat{\alpha}(X_t^{x},\partial_xU(X_t^{x}, \rho_t)))\rmd t.
  \end{equation}
This completes the proof of our desired results.
\qed

At the end of this section, we provide an example of a solvable elliptic master equation 
and demonstrate that its solutions are generally nonunique. Set $b(x,\mu,a)=b_1x+b_2\bar{\mu}+b_3a, f(x,\mu,a)=b_4x\bar{\mu}+Ax^2+Ca^2$,
where $b_1, b_2, b_3, b_4, A, C$ are constants satisfying all conditions in  Theorem \ref{thm: solvable}
and Remark \ref{rm: require detaile}.
Then, we have $\hat{\alpha}(x,y)=-\frac{b_3y}{2C}$ and
\begin{equation}
  H(x,\mu,y)=(b_1x+b_2\bar{\mu})\cdot y+b_4x\bar{\mu}+Ax^2-\frac{b_3^2}{4C}y^2.
\end{equation}
The master equation (\ref{eq: master}) now becomes
\begin{equation}
\label{eq: master linear}
\begin{split}
  r U(x,\mu)=&(b_1x+b_2\bar{\mu})\cdot \partial_{x}U(x,\mu)+b_4x\bar{\mu}+Ax^2\\&-\frac{b_3^2}{4C}(\partial_{x}U(x,\mu))^2
  +\frac{1}{2}\partial_{xx}U(x,\mu)\\&+\tilde{\mathbb{E} }
  \left[\frac{1}{2}\partial_{\tilde{x}}\partial_\mu U(x,\mu,\tilde{\xi })+\partial_\mu U(x,\mu,\tilde{\xi })
  (b_1\tilde{\xi }+b_2\bar{\mu}-\frac{b^2_3}{2C}\partial_{x}U(\tilde{\xi },\mu))  \right].
\end{split}
\end{equation}

We assume that the solution takes the form
\begin{equation}
  U(x,\mu)=a_1x^2+a_2x\bar{\mu}+a_3(\bar{\mu})^2+a_4,
\end{equation}
where $a_1,a_2,a_3,a_4$ are constants that remain to be determined.
Then we have
\begin{equation}
  \begin{split}
    &\partial_xU(x,\mu)=2a_1x+a_2\bar{\mu},\quad\quad\quad \partial_{xx}U(x,\mu)=2a_1,\\
   &\partial_{\mu}U(x,\mu,\tilde{x})= a_2x+2a_3\bar{\mu},\quad\quad \partial_{\tilde{x}}\partial_{\mu}U(x,\mu,\tilde{x})=0.
  \end{split}
\end{equation}
Substituting these into the equation (\ref{eq: master linear}), we obtain
\begin{equation}
  \begin{split}
  &r\left( a_1x^2+a_2x\bar{\mu}+a_3(\bar{\mu})^2+a_4  \right)\\
  =&(b_1x+b_2\bar{\mu})\cdot (2a_1x+a_2\bar{\mu})+b_4x\bar{\mu}+Ax^2-\frac{b_3^2}{4C}(2a_1x+a_2\bar{\mu})^2
  +a_1\\&+ (a_2x+2a_3\bar{\mu})\cdot \left(b_1\bar{\mu}+b_2\bar{\mu}-\frac{b^2_3}{2C} (2a_1\bar{\mu}+a_2\bar{\mu}  )  \right) .
  \end{split}
\end{equation}
Comparing all coefficients, we get the following system of linear and quadratic equations:
\begin{equation}
  \label{eq: roots}
  \begin{cases}
    ra_1=2b_1a_1+A-\frac{b_3^2}{C}a_1^2,\\
    ra_2=2b_2a_1+b_1a_2+b_4-\frac{b_3^2}{C}a_1a_2+a_2(b_1+b_2-\frac{b_3^2}{C}a_1-\frac{b_3^2}{2C}a_2),\\
    ra_3=b_2a_2-\frac{b_3^2}{4C}a_2^2+2a_3(b_1+b_2-\frac{b_3^2}{C}a_1-\frac{b_3^2}{2C}a_2),\\
    ra_4=a_1.
  \end{cases}
\end{equation}
Note that $a_1,a_2,a_3$, and $a_4$
can be solved sequentially, and the equations for $a_1$
and $a_2$
are quadratic.
Since $A>0 $ and $-\frac{b_3^2}{C}<0 $, the equation for $a_1$ must have one positive and one negative root.
If the coefficients are sufficiently well-behaved, for instance, $2b_2a_1+b_4>0$, then the equation for $a_2$
will also have two distinct real roots.
Therefore, the system of equations (\ref{eq: roots}) can have at most four solutions.
It must be emphasized that there exists at most one set of solutions for which the solution of SDE (\ref{eq: master xi state}) 
lies in $L_r^2$, as the $L_r^2$ solution of FBSDE (\ref{eq: mv xi}) is unique.
In a companion article, we establish the uniqueness of solutions to the master equation under certain monotonicity and 
growth assumptions.

\begin{eg}
  We consider a case that $r=2, b_1=b_2=b_4=0, b_3=2, A=2, C=1$, then equation (\ref{eq: master linear}) becomes
\begin{equation}
\begin{aligned}
  2U(x,\mu)=&2x^2-\left(\pa_x U(x,\mu)\right)^2+\frac{1}{2}\pa_{xx}U(x,\mu)\\
  &+\tilde{\mathbb{E} }
  \left[\frac{1}{2}\partial_{\tilde{x}}\partial_\mu U(x,\mu,\tilde{\xi })-2\partial_\mu U(x,\mu,\tilde{\xi })
  \cdot \partial_{x}U(\tilde{\xi },\mu) \right]
\end{aligned}
\end{equation}
Solving the corresponding Equation (\ref{eq: roots}) yields the following solutions:
\begin{equation}
  \begin{aligned}
    &U_1(x,\mu)=\frac{1}{2}x^2+\frac{1}{4},\\
    &U_2(x,\mu)=\frac{1}{2}x^2-3x\bar{\mu}+\frac{3}{2}(\bar{\mu})^2+\frac{1}{4},\\
    &U_3(x,\mu)=-x^2-\frac{1}{2},\\
    &U_4(x,\mu)=-x^2+3x\bar{\mu}-\frac{3}{2}(\bar{\mu})^2-\frac{1}{2}.
  \end{aligned}
\end{equation}
However, to ensure that the solution of Equation (\ref{eq: master xi state}) , 
now written as 
\begin{equation}
  \rmd X_t= -4\pa_xU(X_t,\cL_{X_t}) \rmd t+\rmd B_t,\quad X_0=\xi,
\end{equation}
belong to $L_2^2$, only the solution $U_1(x,\mu)=\frac{1}{2}x^2+\frac{1}{4}$
satisfies this condition.
\end{eg}

\section{Viscosity solution to distribution-dependent elliptic PDE}
\label{sec: viscosity}
\setcounter{equation}{0}

We have proven that the classical solutions of the master equations can be employed to resolve the infinite-time 
FBSDEs. In this section, the solution of the infinite-time FBSDEs will be 
utilized to characterize the viscosity solutions of the distribution-dependent elliptic PDEs.
Specifically, we consider the following FBSDEs with initial state $\xi\in  \mathbb{L} ^2(\mathcal{F} _0)$ and $x\in \mathbb{R} $:
\begin{equation}
  \label{eq: 61}
  \begin{cases}
    \rmd X_{t}^{\xi} = \partial_yH(X_{t}^{\xi},\mathcal{L} _{X_{t}^{\xi}},Y_{t}^{\xi}) \rmd t +\rmd B_{t}, \\
    \rmd Y_{t}^{\xi} = -\left[\partial_x {H} (X_{t}^{\xi},\mathcal{L} _{X_{t}^{\xi}}, Y_{t}^{\xi})-rY_{t}^{\xi}\right] \rmd t + Z_{t}^{\xi}\rmd B_{t}, \\
    X_0^{\xi}=\xi.
    \end{cases}
\end{equation}
\begin{equation}
  \label{eq: 62}
  \begin{cases}
    \rmd X_{t}^{x,\xi} = \partial_yH(X_{t}^{x,\xi},\mathcal{L} _{X_{t}^{\xi}},Y_{t}^{x,\xi}) \rmd t +\rmd B_{t}, \\
    \rmd Y_{t}^{x,\xi} = -\left[\partial_x {H} (X_{t}^{x,\xi},\mathcal{L} _{X_{t}^{\xi}}, Y_{t}^{x,\xi})-rY_{t}^{x,\xi}\right] \rmd t + Z_{t}^{x}\rmd B_{t}, \\
    X_0^x=x.
    \end{cases}
\end{equation}
Here, $H(x,\mu,y)$ is defined in (\ref{eq: hamilton}), and the above equations are the same as (\ref{eq: mv xi1}) and (\ref{eq: mv x1}).

We define $\mcv(x,\mu)\triangleq Y^{x,\xi}_0 $ and 
attempt to prove that it satisfies the following elliptic PDE in the viscosity sense:
\begin{equation}
  \label{eq: pa master}
  \begin{split}
  r \mcu(x,\mu)=&\partial_xH(x,\mu,\mcu(x,\mu))+ \partial_yH(x,\mu,\mcu(x,\mu))\cdot \partial_x\mcu(x,\mu) +
  \frac{1}{2}\partial_{xx}\mcu(x,\mu)\\&+\tilde{\mathbb{E} }
  \left[\frac{1}{2}\partial_{\tilde{x}}\partial_\mu \mcu(x,\mu,\tilde{\xi })+\partial_\mu \mcu(x,\mu,\tilde{\xi })
  \partial_y H(\tilde{\xi },\mu,\mcu(\tilde{\xi},\mu))\right].
\end{split}
\end{equation}
This equation is derived by taking the partial derivative of 
both sides of the master equation (\ref{eq: master}) with respect to $x$. We have the relationship $\mcu(x,\mu)=\pa_ xU(x,\mu)$.

Now, let us give the definition of a viscosity solution for PDE (\ref{eq: pa master}).
We rewrite the PDE as follows:
\begin{equation}
  \label{eq: vis pde}
  (\mathcal{L} \mathcal{U} )[\mathcal{U}](x,\mu)+F(x,\mu,\mathcal{U} (x,\mu))=0,
\end{equation}
where
\begin{equation}
  \begin{split}
  (\mathcal{L}\Phi)[\Psi] (x,\mu)\triangleq&\partial_yH(x,\mu,\Psi(x,\mu))\cdot \partial_x\Phi(x,\mu) +
  \frac{1}{2}\partial_{xx}\Phi(x,\mu)\\&+\tilde{\mathbb{E} }
  \left[\frac{1}{2}\partial_{\tilde{x}}\partial_\mu \Phi(x,\mu,\tilde{\xi })+\partial_\mu \Phi(x,\mu,\tilde{\xi })
  \partial_y H(\tilde{\xi },\mu,\Psi(\tilde{\xi},\mu))\right],
  \end{split}
\end{equation}
and
\begin{equation}
  F(x,\mu,\Psi(x,\mu))\triangleq \partial_xH(x,\mu,\Psi(x,\mu))-r \Psi(x,\mu).
\end{equation}

\begin{defn}
  A function $\Psi\in C^{2,1}(\mathbb{R} \times \mathcal{P} _2)$ 
  is said to be a test function if the quantities:
  \begin{equation}
    \int_\dbR \left\lvert \pa_\mu\Psi(x,\mu,\tilde{x})\right\rvert ^2 \rmd \mu(\tilde{x})
  \end{equation}
  and
  \begin{equation}
    \sup_{\tilde{x}\in \dbR} \left\lvert \pa_{\tilde{x}} \pa_\mu\Psi(x,\mu,\tilde{x})\right\rvert 
  \end{equation}
  are finite, uniformly in $(x,\mu)$ in any compact subset of $\mathbb{R} \times \mathcal{P} _2 $.
\end{defn}

\begin{defn}
  Let $\mcu\in C(\mathbb{R}\times \mathcal{P} _2 )$. Then $\mcu$ is called a viscosity subsolution (resp. supersolution)
  of PDE (\ref{eq: vis pde}) if, whenever $\Psi$ is a test function, and 
  $(x^0,\mu^0)\in \mathbb{R} \times \mathcal{P} _2$ is a local maximum (resp. minimum) of $\mcu-\Psi$, we have
  \begin{equation}
    (\mathcal{L} \Psi )[\cU](x^0,\mu^0)+F(x^0,\mu^0,\mathcal{U} (x^0,\mu^0))\ge0,
  \end{equation}
  (respectively,
  \begin{equation}
    (\mathcal{L} \Psi )[\cU](x^0,\mu^0)+F(x^0,\mu^0,\mathcal{U} (x^0,\mu^0))\le 0,
  \end{equation}
  ).

  The function $\mcu$ is called a viscosity solution of PDE (\ref{eq: vis pde}) if it is 
  both a viscosity subsolution and a viscosity supersolution.
\end{defn}

The proof of $\mcv$ being a viscosity solution to equation (\ref{eq: pa master}) rests on two fundamental results:
\begin{itemize}
  \item The value of $Y^{x,\xi} _0$ depends  solely on the distribution of $\xi$, but not on the specific realization of $\xi$.
  \item $\mcv(x,\mu)\triangleq Y^{x,\xi} _0$ is jointly continuous on $\mathbb{R} \times \mathcal{P} _2$.
\end{itemize}
We will next establish them sequentially in this section and ultimately prove that $\mcv$ is a viscosity solution to equation (\ref{eq: pa master})

\subsection{Yamada-Watanabe theorem for infinite-time FBSDEs}
We shall prove the weak uniqueness for the infinite-time FBSDE:
\begin{equation}
  \label{eq: 6fbsde}
  \begin{cases}
    \rmd X_t=G(t,X_t,Y_t,\mathcal{L}_{X_t})dt+ \rmd B_t, \\
    \rmd Y_t=-F(t,X_t,Y_t,\mathcal{L}_{X_t})dt+Z_t\rmd B_t.
    \end{cases}
\end{equation}
We know that under Assumption \ref{assum: fbsde}, the above FBSDE has a unique solution $(X_t,Y_t,Z_t)$
in $L_r^2$ with initial $X_0\in \dbL^2(\cF_0)$.
We then give the definitions of strong and weak uniqueness.

\begin{defn}[Strong uniqueness]
  We say that the strong uniqueness holds for FBSDE (\ref{eq: 6fbsde}) if on any admissible setup
  $(\Omega,\cF,\dbP,\dbF)$ with input $(X_0, B)$, for any two $\dbF$-progressively measurable $L^2_r$
  three-tuples
  \begin{equation}
    (X^1_t,Y^1_t,Z^1_t)_{t\ge 0},\quad (X^2_t,Y^2_t,Z^2_t)_{t\ge 0}
  \end{equation}
  satisfying the FBSDE (\ref{eq: 6fbsde}) with the same initial condition $X_0$ (up to an exceptional event),
  it holds that
  \begin{equation}
    \dbE\left[\int_0^\infty \rme^{-rt} \left( |X^1_t-X^2_t|^2+
     |Y^1_t-Y^2_t|^2 +|Z^1_t-Z^2_t|^2 \right)\rmd t\right]=0
  \end{equation}
\end{defn}

\begin{defn}[Weak uniqueness]
  For any two set-ups $(\Omega^1, \mathcal{F}^1 ,\mathbb{P}^1,\mathbb{F}^1 )$ and $(\Omega^2, \mathcal{F}^2 ,\mathbb{P}^2,\mathbb{F}^2 )$
  with inputs $(X_0^1,B^1)$ and $(X_0^2,B^2)$, $X_0^1,X_0^2$ have the same law on $\dbR$.
  We say the weak uniqueness holds for FBSDE (\ref{eq: 6fbsde}) if for
  the $L^2_r$ solutions 
  $(X_t^1,Y_t^1,Z^1_t)_{t\ge0}$ and $(X_t^2,Y_t^2,Z^2_t)_{t\ge0}$ on corresponding set-ups,
  the processes
  $(X_t^1,Y_t^1,\int_0^tZ^1_s\rmd s)_{t\ge0}$ and $(X_t^2,Y_t^2,\int_0^tZ^2_s\rmd s)_{t\ge0}$
   have the same distribution.
\end{defn}

We use the same scheme as the one developed by Yamada and Watanabe to
 prove that pathwise uniqueness of solutions of FBSDE's implies uniqueness in the sense of
 probability law.

\begin{thm}
  Assume that on $(\Omega, \mathcal{F} ,\mathbb{P},\mathbb{F} )$ with input $(X_0, B_t)_{t\ge 0}$, the FBSDE
  has a unique strong solution $(X_t,Y_t,Z_t)$.
  Then, the law of $(X_0,B_t,X_t,Y_t,\int_0^tZ_s\rmd s)$ only depends on $\mathcal{L} (X_0)$.
\end{thm}

\proof 
First, we denote by $C([0,\infty);\dbR)$ the space of continuous $\dbR$-valued functions on $[0,\infty)$
equipped with the metric of uniform convergence on compacts:
\begin{equation}
  d(\omega^1,\omega^2)=\sum_{n\ge0}2^{-n}\sup_{t\in[0,n]}\max(|\omega_t^1-\omega_t^2|,1).
\end{equation}
And define
\begin{equation}
  \begin{split}
  &\Oin\triangleq\dbR\times C([0,\infty);\dbR),\\ &\Oout\triangleq C([0,\infty);\dbR)\times C([0,\infty);\dbR)\times C([0,\infty);\dbR),\\
&\Ocan\triangleq \Oin\times \Oout.
  \end{split}
\end{equation}

Let us consider  two filtered probability spaces $(\Omega^i, \mathcal{F}^i ,\mathbb{P}^i,\mathbb{F}^i )$
with identically distributed inputs $(X_0^i,B^i)$, $i=1,2$,
on each of which a solution $(X_t^i,Y_t^i,\int_0^tZ^i_s\rmd s)_{t\ge0}$
to the FBSDE (\ref{eq: 6fbsde}) is defined.
Denote by $Q^1$ and $Q^2$ the distribution of $(X_0^1, B_t^1,X_t^1,Y_t^1,\int_0^tZ^1_s\rmd s)_{t\ge0} $
and $(X_0^2, B_t^2,X_t^2,Y_t^2,\int_0^tZ^2_s\rmd s)_{t\ge0} $ on $\Ocan=\Oin\times \Oout$,
by $\Qin$ the common distribution of the processes $(X_0^1,B_t^1), (X_0^2, B_t^2)$ on $\Oin$.

Let us now define for $i\in\{1,2\}$,
\begin{equation}
  q^i(\oin;F):\Oin\times \mathcal{B} (\Oout)\rightarrow [0,1]
\end{equation}
 as the regular conditional probability for $\mathcal{B} (\Oout)$ given $\oin\in \Oin$ (under $Q^i$). It satisfies:
\begin{itemize}
  \item $\forall \oin\in \Oin$, $q^i(\oin; \cdot)$ is a probability measure on $(\Oout,\mathcal{B} (\Oout))$.
  \item $\forall F\in \mathcal{B} (\Oout)$, the mapping $(\xi,w)\rightarrow q^i(\xi,w;F)$ is $\cB(\dbR)\otimes \cB(C([0,\infty),\dbR))$-measurable.
  \item $\forall F\in \cB(\Oout), \forall G\in \cB(\Oin)$:
        \begin{equation}
          Q^i(G\times F)=\int _G q^i(\oin;F)\Qin(\rmd \oin).
        \end{equation}
\end{itemize}

Next, we need an enlarged space $(\Oto, \mathcal{G}, Q )$ to support all processes. Define
\begin{equation}
  \Oto\triangleq \Oin \times\Oout\times\Oout,
\end{equation}
 and $\mathcal{G} $ is the
 completion of the $\sigma$-field $\cB(\Ocan)\otimes \cB(\Oout)$ by the collection $\cN$ of all null sets under the
 probability measure
 \begin{equation}
  Q(G\times F_1\times F_2)=\int_G q^1(\oin;F_1)q^2(\oin;F_2)\Qin(\rmd \oin),
 \end{equation}
where $F_1,F_2\in\cB(\Oout), G\in \cB(\Oin)$.

We observe that $Q(G\times F_1\times\Oout)=Q^1(G\times F_1), Q(G\times \Oout\times F_2)=Q^2(G\times F_2).$
Inparticular, we denote by $(\xi,w,x^1,y^1,\zeta^1,x^2,y^2,\zeta^2)$ the canonical process on $\Oto$,
then $(\xi,w,x^1,y^1,\zeta^1)$ has distribution $Q^1$ and $(\xi,w,x^2,y^2,\zeta^2)$ has distribution $Q^2$.

We define $(z^i_t)_{t\ge0},i\in\{1,2\}$ by
\begin{equation}
  z^i_t=\begin{cases}
 \lim_{n\rightarrow \infty} n\left( \zeta^i_t-\zeta^i_{(t-\frac{1}{n})+}  \right)\quad  ~\mbox{if the limit exists,}
  \\ 0\quad ~\mbox{otherwise.}
   \end{cases}
\end{equation}
Since $Q^i$-a.s., $\int_0^t Z^i_s\rmd s$ is absolutely continuous, we know that $Q$-a.s.,$\zeta^i_t$ is absolutely continuous.
So
\begin{equation}
  \zeta^i_t=\int_0^tz^i_s\rmd s,\quad t\ge0.
\end{equation}
Moreover, we have
\begin{equation}
  \begin{split}
  \dbE^Q\int_0^\infty \rme^{-rt} |z^i_t|\rmd t 
  &\le  \liminf_{n\rightarrow \infty}\dbE^Q\int_0^\infty \rme^{-rt}|n(\zeta^i_t-\zeta^i_{(t-\frac{1}{n})+})|^2\rmd t\\
  &=  \liminf_{n\rightarrow \infty}n^2\dbE^{\dbP^i}\int_0^\infty \rme^{-rt}|\int_{(t-\frac{1}{n})+}^tZ_s^i\rmd s |^2\rmd t\\
  &\le \liminf_{n\rightarrow \infty}n\dbE^{\dbP^i}\int_0^\infty \rme^{-rt}\int_{(t-\frac{1}{n})+}^t|Z_s^i|^2\rmd s \rmd t\\
  &\le \dbE^{Q^i}\int_0^\infty \rme^{-rs}|Z_s^i|^2\rmd s,
  \end{split}
\end{equation}
the last inequality following from  Fubini's theorem.

Let us now endow $(\Oto, \cG, Q)$ with the filtration $\dbG$, where
$\dbG=\{\cG_t\}_{t\ge 0}$ is the complete and right-continuous augmentation under $Q$
of the canonical filtration
\begin{equation}
  \cH_t=\{(\xi,w_s,x^1_s,y^1_s,\zeta^1_s,x^2_s,y^2_s,\zeta^2_s);0\le s\le t\}
\end{equation}
on $\Oto$.
It is easy to see that $\xi$ is $\cG_0$-measurable and that $(w_s,x^1_s,y^1_s,z^1_s,x^2_s,y^2_s,z^2_s)_{t\ge 0}$
is $\{\cG_t\}_{t\ge 0}$-progressively measurable. Moreover, for $i\in\{1,2\}$,
\begin{equation}
  Q\left\{\omega\in \Oto; (x,w,x^i,y^i,\zeta^i)\in A\right\}
  =Q^i\left\{ (X_0^i,B^i,X^i,Y^i,\int_0^\cdot Z^i_s\rmd s)\in A\right\};\quad A\in\cB(\Ocan).
\end{equation}

Actually, we just have to prove that $(w_t)_{t\ge 0}$ is a $\{\cG_t\}_{t\ge 0 }$ Brownian motion:
we follow the proof given in \cite{DELARUE2002209} (Remark 1.6).
Let us firstly define
\begin{equation}
  \pi_t: C([0,\infty);\dbR)\rightarrow C([0,t];\dbR),\quad  h\rightarrow h|_{[0,t]},
\end{equation}
and 
\begin{equation}
  \pi'_t:\gamma\rightarrow \gamma_t,\quad h\rightarrow h|_{[0,t]},
\end{equation}
where $\gamma=\Oout$ and 
\begin{equation}
  \gamma_t=C([0,t];\dbR)\times C([0,t];\dbR)\times C([0,t];\dbR).
\end{equation}
Endowing $C([0,t];\dbR)$ and $\gamma_t$ with their borelian $\sigma$-fields, we define
\begin{equation}
  \cK_t\triangleq \sigma\{\pi_t\},\quad \cK'_t\triangleq \sigma\{\pi'_t\}.
\end{equation}
Using the separability of the spaces $C([0,t];\dbR)$ and $\gamma_t$,
we see that
\begin{equation}
  \cK_t=\sigma\{w_s; 0\le s \le t\},
\end{equation}
and that $\forall i\in\{1,2\}, \forall A\in \cK'_t$, the set $\{(X^i,Y^i,\int_0^\cdot Z^i_s\rmd s)\in A \}$ belongs to $\cF^i_t$.

 Now, considering $A\in\cK'_t$ , we want to show that, for $i\in\{1,2\}$, the map
 \begin{equation}
  \Oin\rightarrow [0,1],\quad (\xi,w)\rightarrow q^i(\xi,w;A)
 \end{equation}
is measurable with respect to the completion of the $\sigma$-field $\cB(\dbR)\otimes \cK_t$ under the
probability measure $\Qin$, denoted $\overline{\cB(\dbR)\otimes \cK_t }$.
Indeed, let us consider $F\in\cB(\dbR)$, $G_1\in\cK_t$  and $G_2\in \sigma \{w_s-w_t; s\ge t\}$. Then, $\forall i\in\{1,2\}$
\begin{equation}
\begin{split}
  &\int I_F(\xi)I_{G_1}(w)I_{G_2}(w)q^i(\xi,w;A)\Qin(\rmd \xi \rmd w)\\
=& \dbE^{\dbP^i}\left[I_F(X_0)I_{G_1}(B^i)I_{G_2}(B^i)I_A(X^i,Y^i,\int_0^{\cdot}Z_s^i\rmd s)   \right]\\
=& \dbE^{\dbP^i}\left[I_F(X_0)I_{G_1}(B^i)I_A(X^i,Y^i,\int_0^{\cdot}Z_s^i\rmd s)   \right]
     \dbE^{\dbP^i}\left[I_{G_2}(B^i)  \right]\\
=&\int I_F(\xi)I_{G_1}(w)q^i(\xi,w;A)\Qin(\rmd \xi \rmd w)\int I_{G_2}(w)\Qin(\rmd \xi \rmd w).
\end{split}
\end{equation}
Hence, the map $(\xi,w)\rightarrow q^i(\xi,w;A)$ is measurable with respect to $\overline{\cB(\dbR)\otimes \cK_t }$.

Now we prove that $(w_t)_{t\ge 0}$ is a $\{\cG_t\}_{t\ge 0 }$ Brownian motion.
Let us consider $(A,A')\in(\cK'_t)^2, F\in \cB(\dbR), G_1\in\cK_t$  and $G_2\in \sigma \{w_s-w_t; s\ge t\}$.
Then,
\begin{equation}
\begin{split}
  &\dbE^Q \left[ I_F(\xi)I_{G_1}(w)I_{G_2}(w)I_A(x^1,y^1,\zeta^1)I_{A'}(x^2,y^2,\zeta^2)   \right]\\
=&\int I_F(\xi)I_{G_1}(w)I_{G_2}(w)q^1(\xi,w;A)q^2(\xi,w;A')\Qin(\rmd \xi \rmd w)\\
=&\int I_F(\xi)I_{G_1}(w)q^1(\xi,w;A)q^2(\xi,w;A')\Qin(\rmd \xi \rmd w)\int I_{G_2}(w)\Qin(\rmd \xi \rmd w)\\
=&\dbE^Q \left[ I_F(\xi)I_{G_1}(w)I_A(x^1,y^1,\zeta^1)I_{A'}(x^2,y^2,\zeta^2)   \right]
   \dbE^Q \left[I_{G_2}(w)\right].
\end{split}
\end{equation}
Noting that $\cH^t=\cB(\dbR)\otimes \cK_t\otimes \cK'_t\otimes \cK'_t$, we conclude that $(w_t)_{t\ge 0}$ is a $\{\cG_t\}_{t\ge 0 }$ Brownian motion.

At last, applying the same procedure to $(z^i_t)_{t\ge 0}, i\in\{1,2\}$ in \cite{carmona2018probabilistic} (Volume II, Lemma 1.27), we obtain that $Q-a.s.$: for all $0\le t\le T,i\in\{1,2\}$,
\begin{equation}
  \begin{cases}
    x^i_t=\xi+\int_0^t G(s,x^i_s,y^i_s,\mathcal{L}_{x^i_s})ds+  w_t, \\
   y^i_t=y^i_T+\int_t^T F(s,x^i_s,y^i_s,\mathcal{L}_{x^i_s})ds-\int_t^T z^i_s\rmd w_s,\\
    \dbE^Q\left[\int_0^\infty \rme^{-rt} \left( |x^i_t|^2+
     |y^i_t|^2 +|z^i_t|^2 \right)\rmd t\right]<\infty.
  \end{cases}
\end{equation}
Through the strong uniqueness, we know that under $Q$, 
the processes $(x^1_t,y^1_t,\zeta^1_t)_{t\ge 0}$ and $(x^2_t,y^2_t,\zeta^2_t)_{t\ge 0}$
have the same law, then we get the desired result.

\qed

\begin{rem}
  For any two initial states  $\xi_1,\xi_2\in  \mathbb{L} ^2(\mathcal{F} _0)$ with the same distribution and $x\in\mathbb{R} $, 
  the
solutions $Y^{x,\xi_1}, Y^{x,\xi_2}$ to (\ref{eq: 62}) are the same.
We can say that  $Y^{x,\xi} _0$ depends  solely on the distribution of $\xi$.
\end{rem}

\subsection{Connection with distribution-dependent elliptic PDE}

For FBSDEs (\ref{eq: 61}) and (\ref{eq: 62}), we define $\mcv(x,\mu)\triangleq Y^{x,\xi}_0 $ for some initial state
  $\xi$ with distribution $\mu$.  
We first present the continuity result for $\mcv$.

\begin{lem}
Assuming all assumptions in Theorem \ref{thm: solvable} are satisfied,
   we have, for any $x,x'\in\mathbb{R} $ and $\mu,\mu'\in \mathcal{P} _2$, there exists a constant $C>0$ such that
  \begin{equation}
    |\mcv(x,\mu)-\mcv(x',\mu')\le C(|x-x'|+\mathcal{W} _2(\mu,\mu'))
  \end{equation}
\end{lem}
\proof 
Let $(X^{\xi_1},Y^{\xi_1},X^{x_1,\xi_1},Y^{x_1,\xi_1})$ and $(X^{\xi_2},Y^{\xi_2},X^{x_2,\xi_2},Y^{x_2,\xi_2})$
be the solutions of Equations (\ref{eq: 61}) and (\ref{eq: 62}) with $\mathcal{L} _{\xi_1}=\mu_1,\mathcal{L} _{\xi_2}=\mu_2$.
Set 
\begin{equation}
  \begin{split}
    &\hat{X}^{\xi}=X^{\xi_1}-X^{\xi_2},\quad \hat{Y}^{\xi}=Y^{\xi_1}-Y^{\xi_2},\\
    &\hat{X}^{x,\xi}=X^{x_1,\xi_1}-X^{x_2\xi_2},\quad \hat{Y}^{x,\xi}=Y^{x_1,\xi_1}-Y^{x_2,\xi_2}.
  \end{split}
\end{equation}
$C_1,C_2,C_3,C_4$ appeared in the following proof are positive constants.

Applying It\^{o}'s formula to $\rme^{-rt}|Y^{x_1,\xi_1}_t-Y^{x_2,\xi_2}_t|^2$, we get
\begin{equation}
  \begin{split}
    |Y^{x_1,\xi_1}_0-Y^{x_2,\xi_2}_0|^2&\le 
    C_1\mathbb{E} \int_0^{\infty} \rme^{-rt}\left[(\hat{X}^{x,\xi}_t)^2+ (\hat{Y}^{x,\xi}_t)^2+W_2^2(X^{\xi_1}_t,X^{\xi_2}_t)  \right]\rmd t\\
    &\le C_1\mathbb{E}\int_0^{\infty} \rme^{-rt}\left[(\hat{X}^{x,\xi}_t)^2+ (\hat{Y}^{x,\xi}_t)^2+(\hat{X}^{\xi}_t )^2  \right]\rmd t.
  \end{split}
\end{equation}
Applying It\^{o}'s formula to $\rme^{-rt}\hat{X}^{x,\xi}\hat{Y}^{x,\xi}$, we get
\begin{equation}
  \begin{split}
    \mathbb{E}\int_0^{\infty} \rme^{-rt}\left[(\hat{X}^{x,\xi}_t)^2+ (\hat{Y}^{x,\xi}_t)^2 \right]\rmd t
    \le C_2\left( \hat{X}^{x,\xi}_0\hat{Y}^{x,\xi}_0 + \mathbb{E}\int_0^{\infty} \rme^{-rt}(\hat{X}^{\xi}_t )^2\rmd t  \right).
  \end{split}
\end{equation}
It turns out that
\begin{equation}
  \begin{split}
  (\hat{Y}^{x,\xi}_0)^2\le& C_1C_2\hat{X}^{x,\xi}_0\hat{Y}^{x,\xi}_0+(C_1C_2+C_1)\mathbb{E}\int_0^{\infty} \rme^{-rt}(\hat{X}^{\xi}_t )^2\rmd t\\
    \le&\frac{1}{2}(\hat{Y}^{x,\xi}_0)^2+\frac{C_1^2C_2^2}{2}(\hat{X}^{x,\xi}_0)^2+(C_1C_2+C_1)\mathbb{E}\int_0^{\infty} \rme^{-rt}(\hat{X}^{\xi}_t )^2\rmd t.
  \end{split}
\end{equation}
Thus,
\begin{equation}
  (\hat{Y}^{x,\xi}_0)^2\le C_1^2C_2^2(\hat{X}^{x,\xi}_0)^2+2(C_1C_2+C_1)\mathbb{E}\int_0^{\infty} \rme^{-rt}(\hat{X}^{\xi}_t )^2\rmd t.
\end{equation}
Applying the same arguments to $\rme^{-rt}|\hat{Y}^{\xi}_t|^2$ and $\rme^{-rt}\hat{X}^{\xi}\hat{Y}^{\xi}$, we can get that
\begin{equation}
  \mathbb{E}\int_0^{\infty} \rme^{-rt}(\hat{X}^{\xi}_t )^2\rmd t\le C_3 \mathbb{E}|\xi_1-\xi_2|^2.
\end{equation}
So, we can deduce that 
\begin{equation}
  |Y^{x_1,\xi_1}_0-Y^{x_2,\xi_2}_0|^2\le C_4\left(|x_1-x_2|^2+\mathbb{E}|\xi_1-\xi_2|^2\right).
\end{equation}
Since $Y_0^{x,\xi}$ depends solely on the distribution of $\xi$, we have
\begin{equation}
  |\mcv(x_1,\mu_1)-\mcv(x_2,\mu_2)|^2\le C_4\left(|x_1-x_2|^2+\mathcal{W} _2^2(\mu_1,\mu_2)\right).
\end{equation}
Then, we get the desired result.
\qed

We now assert the following theorem.
\begin{thm}
  Considering all assumptions in Theorem \ref{thm: solvable} are satisfied,
  the function $\mcv(x,\mu)$ is a viscosity solution of PDE (\ref{eq: vis pde}).
\end{thm}
\proof 
We only show that $\mcv$ is a viscosity subsolution of PDE (\ref{eq: vis pde}). 
A similar argument will show that it is also a viscosity supersolution of PDE (\ref{eq: vis pde}).

Due to the uniqueness of solutions of FBSDEs, it is not hard to see that for any 
$t\ge0$, $\mcv(X_t^{x,\xi},\mathcal{L} _{X_t^{\xi}})=Y_t^{x,\xi}$. Let $\Psi$ be a test function
and  $(x^0,\mu^0)\in \mathbb{R} \times \mathcal{P} _2$ be a local maximum  of $\mcu-\Psi$. 
We assume without loss of generality that $\mcv(x^0,\mu^0)=\Psi(x^0,\mu^0)$.
We suppose that
\begin{equation}
  (\mathcal{L} \Psi )[\mcv](x^0,\mu^0)+F(x^0,\mu^0,\mathcal{V} (x^0,\mu^0))<0,
\end{equation}

It follows from the above that there exists an open subset $O\subset \mathbb{R} \times \mathcal{P} _2$ that contains
$(x^0,\mu^0)$ such that for all $(x,\mu)\in O$,
\begin{equation}
  \begin{cases}
    \mcv(x,\mu)\le \Psi(x,\mu),\\
    (\mathcal{L} \Psi )[\mcv](x,\mu)+F(x,\mu,\mathcal{V} (x,\mu))<0.
  \end{cases}
\end{equation}
Taking an initial state $\xi^0\in\mathbb{L} ^2(\mathcal{F} _0;\mu^0)$, we consider the processes 
$(X_t^{\xi^0},Y_t^{\xi^0},Z_t^{\xi^0})$ and $(X_t^{x^0,\xi^0},Y_t^{x^0,\xi^0},Z_t^{x^0,\xi^0})$, which are solutions to FBSDEs (\ref{eq: 61}) and (\ref{eq: 62}), respectively.
For some $T>0$, let $\tau$ denote the stopping time
\begin{equation}
  \tau\triangleq \inf\{t>0|(X_t^{x^0,\xi^0},\mathcal{L} _{X_t^{\xi^0}})\notin O\}\land T.
\end{equation}

Let $\rho_t\triangleq \cL_{X_t^{\xi^0}}$,
we first note that the pair of processes
\begin{equation}
  (
\bar{Y}\left(t\right),\bar{Z}\left(t\right)
  )\triangleq
(
Y^{x^0,\xi^0}_{t\wedge\tau},I_{\left[0,\tau\right]}(t)Z^{x^0,\xi^0}_t
),\quad0\leq t\leq T,
\end{equation}
is the solution of the BSDE
\begin{equation}
  \begin{split}
    \bar{Y}_t=&\mcv(X^{x^0,\xi^0}_{\tau},\rho_{\tau})
    +\int_t^TI_{[0,\tau]}(s)F(X^{x^0,\xi^0}_{s}, \rho_s, \mcv(X^{x^0,\xi^0}_{s}, \rho_s) )\rmd s
    \\&-\int_t^T\bar{Z}_s\rmd B_s, \quad\quad 0\le t\le T.
  \end{split}
\end{equation}
Next, it follows from It\^{o}'s formula that the pair of processes
\begin{equation}
(\hat{Y}_t,\hat{Z}_t)\triangleq (\Psi(X^{x^0,\xi^0}_{t\wedge \tau},\rho_{t\wedge \tau}  ), 
I_{\left[0,\tau\right]}(t)\partial_x \Psi( X^{x^0,\xi^0}_{t\wedge \tau},\rho_{t\wedge \tau}   )  ) 
\end{equation}
is the solution of the BSDE
\begin{equation}
    \begin{split}
    \hat{Y}_t=&\Psi(X^{x^0,\xi^0}_{\tau},\rho_{\tau})
    -\int_t^TI_{[0,\tau]}(s) \mathcal{L} \Psi[\mcv] (X^{x^0,\xi^0}_{s}, \rho_s )\rmd s
    \\&-\int_t^T\hat{Z}_s\rmd B_s, \quad\quad 0\le t\le T.
  \end{split}
\end{equation}
Define
\begin{equation}
  \beta_s=-\mathcal{L} \Psi[\mcv] (X^{x^0,\xi^0}_{s}, \rho_s )
  -F(X^{x^0,\xi^0}_{s}, \rho_s, \mcv(X^{x^0,\xi^0}_{s}, \rho_s) )
\end{equation}
and 
\begin{equation}
  \left(\tilde{Y}\left(t\right),\tilde{Z}\left(t\right)\right)=\left(\hat{Y}\left(t\right)-\bar{Y}\left(t\right),\hat{Z}\left(t\right)-\bar{Z}\left(t\right)\right).
\end{equation}
We have
\begin{equation}
  \begin{split}
  \tilde{Y}_t=&\Psi(X^{x^0,\xi^0}_{\tau},\rho_{\tau})-\mcv(X^{x^0,\xi^0}_{\tau},\rho_{\tau})
  +\int_t^TI_{[0,\tau]}(s)\beta_s\rmd s\\&-\int_t^T\tilde{Z}_s\rmd B_s,\quad \quad 0\le t\le T.
  \end{split}
\end{equation}
Therefore,
\begin{equation}
  \tilde{Y}_0=\mathbb{E} \left[\tilde{Y}_{\tau}+\int_0^{\tau}\beta_s \rmd s  \right]
\end{equation}
Now, from the choice of $O$ and $\tau$, a.s.
\begin{equation}
  \tilde{Y}_{\tau}=\Psi(X^{x^0,\xi^0}_{\tau},\rho_{\tau})-\mcv(X^{x^0,\xi^0}_{\tau},\rho_{\tau})\ge 0,
  \quad\beta_s>0, \quad s\in [0,\tau].
\end{equation}
Consequently, $ \tilde{Y}_0=\Psi(x^0,\mu^0)-\mcv(x^0,\mu^0)>0$, which contradicts the earlier assumption.
\qed

\section*{Acknowledgements}
Song Y. is financially supported by National Key R\&D Program of China (No. 2024YFA1013503 \& No. 2020YFA0712700) and the National Natural Science Foundation of China (No. 12431017).

\end{document}